\documentclass[12pt,leqno]{amsart}
\usepackage{amsmath,amssymb,txfonts}
\usepackage{amssymb}
\usepackage{amsxtra}
\usepackage{amsmath}
\usepackage{txfonts}

 \usepackage{color}

\usepackage[cp1252]{inputenc}

\usepackage{mathrsfs}
\usepackage{graphicx}

 \usepackage[active]{srcltx}

\allowdisplaybreaks

\def\rr{{\mathbb R}}
\def\rn{{{\rr}^n}}

\def\fz{\infty}

\def\dist{{\mathop\mathrm{\,dist\,}}}
\def\loc{{\mathop\mathrm{\,loc\,}}}

\def\bz{\beta}

\def\diam{{\mathop\mathrm{\,diam\,}}}

\def\r{\right}
\def\lf{\left}

\newtheorem{thm}{Theorem}[section]
\newtheorem{lem}{Lemma}[section]

\newtheorem{rem}{Remark}[section]

\numberwithin{equation}{section}

\begin{document}
\arraycolsep=1pt

\title[A $(\phi_\frac{n}{s}, \phi)$-Poincar\'e inequality in John domain ]{
A $(\phi_\frac{n}{s}, \phi)$-Poincar\'e inequality in John domain}

\author{Shangying Feng}
\author{Tian Liang}
\address{ Shangying Feng: Department of Mathematics, Beijing Normal University, Beijing 100191, P.R. China;}
\email{202021130044@mail.bnu.edu.cn}
\address{ Tian Liang: School of Mathematics and Statisties, Huizhou University, Guangdong 516007, P.R. China;}
\email{liangtian@hzu.edu.cn}
\thanks{     }

\date{\today }
\maketitle

\begin{center}
\begin{minipage}{13.5cm}\small{\noindent{\bf Abstract}\quad
Let  $\Omega$ be a bounded domain in $\rn$ with $n\ge2$ and $s\in(0,1)$.
Assume that  $\phi : [0, \fz) \to [0, \infty)$  be  a  Young function obeying the doubling condition with  the   constant  $K_\phi<2^{\frac{n}{s}}$.
We demonstrate that $\Omega $ supports a
$(\phi_\frac{n}{s}, \phi)$-Poincar\'e   inequality
if it is  is a John domain.
 Alternately,  assume further that $\Omega$ is a bounded domain that is quasiconformally equivalent to some uniform domain when  $n\ge3$ or  a simply connected domain when $n=2$.
We demonstrate $\Omega$ is a John domain if a $(\phi_\frac{n}{s}, \phi)$-Poincar\'e  inequality holds.
}
\end{minipage}
\end{center}

\section{Introduction\label{s1}}
Let   $n\geq2$ and $\Omega \subset \mathbb{R}^n $ be a bounded
 domain.
Suppose that $\phi$ is a Young function in $[0, \infty)$,  that is,  $\phi \in C[0, \infty)$ is convex and satisfies $\phi(0)=0,  \phi(t)>0$ for $t>0$ and $\lim_{t \to \infty}\phi(t)=\infty$. For any s $\in (0, 1)$,  define the intrinsic fractional Orlicz-Sobolev space $\dot{V}_\ast^{s, \phi}(\Omega)$ as the collection of all measurable functions $u$ in $\Omega$ with the semi-norm
\begin{eqnarray*}
\|u\|_{{\dot {V}_\ast}^{s, \phi}(\Omega)} := \inf \left\{\lambda > 0 : \int_{\Omega} \int_{|x-y|<\frac{1}{2}d(x, \partial\Omega)}  \phi \left(\frac{|u(x)-u(y)|}{\lambda|x-y|^{s}}\right) \frac{dxdy}{|x-y|^{n}}\leq 1 \right\}<\infty.
\end{eqnarray*}
Modulo constant functions,  $\dot{V}_\ast^{s, \phi}(\Omega)$ is a Banach space.
When $s=1$, we usually consider the classical Orlicz-Sobolev space $W^{1, \phi}(\Omega)$,
whose  sharp embedding  has been solved in \cite{c96}(see also \cite{c97} for
an alternate formulation of the solution).

Alberico et al. \cite{refc3} established  an imbedding of $\dot{V}^{s, \phi}_\ast(\rn) $
into certain  Orlicz target space.
Recall that for any Young function $\psi$,  the Orlicz space   $L^\psi(\Omega)$  is the collection of all   $u\in L^1_\loc(\Omega)$ whose norm
\begin{equation*}
\|u\|_{L^{\psi}(\Omega)}:=\inf \left\{\lambda > 0 : \int_{\Omega}  \psi \left(\frac{|u|}{\lambda}\right) dx\leq 1 \right\}<\infty.
\end{equation*}
The following is a more thorough description.
\begin{thm}  Let  $\phi$ be a Young function  satisfying
\begin{equation}\label{eq2}
\int_0^t(\frac{\tau}{\phi(\tau)})^\frac{s}{n-s}d\tau<\infty
\end{equation}
and
\begin{equation}\label{eq3}
\int_0^\infty(\frac{\tau}{\phi(\tau)})^\frac{s}{n-s}d\tau=\infty.
\end{equation}
Define    $\phi_\frac{n}{s}:=\phi\circ H^{-1}$,
where
\begin{equation}\label{eq5}
H(t)=\left(\int_0^t(\frac{\tau}{\phi(\tau)})^\frac{s}{n-s}d\tau\right)^\frac{n-s}{n}~\forall ~t\geq0.
\end{equation}
Then $V_\ast^{s,\phi}(\rn)\subset L^{\phi_{n/s}}(\rn)$, that is,
for any $u\in V_\ast^{s,\phi}(\rn)$ with $\left|\{x\in\rn||u(x)|>t\}\right|<\fz$ for every $t>0$, one has $u\in  L^{\phi_{n/s}}(\rn)$ with
$\|u\|_{L^{\phi_{n/s}}(\rn)}\le C \|u\|_{V_\ast^{s,\phi}(\rn)}$, where $C$ is a constant independent of $u$.
\end{thm}
They also showed that $L^{\phi_{n/s}} (\rn)$ is optimal target spaces for the imbeding of
$\dot{V}_{*}^{s, \phi}(\rn) $  in the sense that
  if    $\dot{V}_{*}^{s, \phi}(\rn) \subset L^{A}(\rn)$ holds  for another Orlicz space $L^{A}(\rn)$,  then $L^{\phi_{n/s}}(\rn) \subset L^{A}(\rn)$.


%

We are interested in bounded domains which supports
the imbedding $V_\ast^{s,\phi}(\Omega)\subset L^{\phi_{n/s}}(\Omega)$ or  $(\phi_\frac{n}{s}, \phi)$-Poincar\'e inequality,  that is,
 there exists a constant $C\geq1$ such that
\begin{align}\label{ueq1}
\|u-u_\Omega\|_{L^{\phi_\frac{n}{s}}(\Omega)}\leq C\|u\|_{{\dot {V}_\ast}^{s, \phi}(\Omega)},
\end{align}
for every $u \in L^1(\Omega)$,
where $u_E =  \fint_{E}u = \frac{1}{|E|} \int _{E} u dx $ denotes the average of $u$
in the set of $E$ with $\lvert E\rvert>0$.

The major finding of this article is to characterize
the imbedding   $V_\ast^{s,\phi}(\Omega)\subset L^{\phi_{n/s}}(\Omega)$
via John domains under specific doubling assumption in $\phi$;
see   Theorem 1.2 below.
Remember that a bounded domain $\Omega \subset \rn$ is called as a $c$-John domain with respect to some $x_0 \in \Omega $ for some $c>0$ if for  each $x \in\Omega$,  there is a rectifiable curve $\gamma : [0, T] \rightarrow \Omega$ parameterized by arc-length such that $\gamma(0)=x$,  $\gamma(T)=x_0$ and $d(\gamma(t),  \Omega ^{\complement}) > ct $ for all $t>0$.
For further research on $c$-John domains, see \cite{m79, r83, bsk96, m88, b82, bk95, bkl95}
and references therein.
We say that a Young function $\phi$ has the doubling property ($\phi \in \Delta_2$) if
\begin{equation}\label{da2}
K_{\phi}:=\sup_{t>0}\frac{\phi(2t)}{\phi(t)}<\infty.
\end{equation}
Note that if a Young function $\phi \in \Delta_2$ with $K_\phi<2^{\frac{n}{s}}$,  then $\phi$ satisfies \eqref{eq2}  and \eqref{eq3};  see Lemma \ref{le1}.

\begin{thm}\label{th1}
 Let $0 < s < 1$. Suppose $\phi$ is a Young function and $\phi \in \Delta_2$ with $K_\phi<2^{\frac{n}{s}}$ in \eqref{da2}.
\begin{enumerate}
\item[(i)] If $\Omega \subset \mathbb{R}^n $ is a $c$-John domain,
then $\Omega$ supports the $(\phi_\frac{n}{s}, \phi)$-Poincar\'e inequality \eqref{ueq1} with the constant $C$ depending on $n,  s,  c$ and $K_\phi$.

\item[(ii)] Assume further that $\Omega \subset \mathbb{R}^n $ is a bounded simply connected planar domain,  or a bounded domain which is a quasiconformally equivalent to some uniform domain when $n\geq3$. If $\Omega$ supports the $(\phi_\frac{n}{s}, \phi)$- Poincar\'e inequality,  then $\Omega$ is a c-John domain,  where the constant $c$ depend on $n,  s,  C, K_\phi$ and $\Omega$.
\end{enumerate}

\end{thm}

Theorem 1.2 extends several known results in the literature; for details see the following remark.

\begin{rem}\rm

(i) For $1\le p<n$, $c$-John domain  $\Omega$ supports Sobolev $\dot{W}^{1,p}$-imbedding or $(\frac{np}{n-p},p)$-Poincar\'e inequality:
\begin{align}\label{1p}
\|u- u_{\Omega}\|_{L^{np/(n-p)}(\Omega)} \leq C \|u\|_{\dot{ W }^{1, p}(\Omega)}  \quad \forall u \in \dot{ W }^{1, p}(\Omega),
\end{align}
where the constant $C$ depends on $n,p$ and $c$;
see Reshetnyak \cite{r83} and Martio \cite{m88} for $1<p<n$ and
Borjarski \cite{b89} (and also Hajlasz \cite{h01}) for $p=1$.
Conversely,  further assume that $\Omega$ is a
bounded simply connected planar domain or a domain that is quasiconformally equaivalently to
some uniform domain when $n\ge3$.
Buckley and Koskela \cite{bk95} proved that   if   \eqref{1p} holds,
then $\Omega$ is a $c$-John domain.

(ii) For $0<s<1$ and $1\le p<\fz$,
the intrinsic fractional Sobolev space $\dot{W}_\ast^{s, p}(\Omega) $  consists of
 all functions  $u\in L^1_\loc(\Omega)$ with
  \begin{eqnarray*}
\|u\|_{\dot{ W }_\ast^{s, p}(\Omega)} :=   \lf(\int_{\Omega} \int_{|x-y|<\frac{1}{2}d(x, \partial\Omega)}  \frac{|u(x)-u(y)|^p}{|x-y|^{n+sp}} \, dxdy\r)^{1/p}<\fz.
\end{eqnarray*}
In the special case  $ \phi(t)=t^p$ with $p\ge1 $,
  $\dot{V}_\ast^{s, \phi}(\Omega)$ is exactly $ \dot{ W }_\ast^{s, p}(\Omega)$.
%

For   $s\in(0, 1)$ and $1\le p<n/s$, 
Dyda-Ihnatsyeva-V\"{a}h\"{a}kangas \cite{d16} for $p=1$ and 
Hurri-Syrj\"{a}nen-V\"{a}h\"{a}kangas \cite{h13} for $1<p<n/s$
proved  that $c$-John domain $\Omega$ supports the  following fractional
  $(np/(n-sp),  p)_{s}$-Poincar\'e   inequality (or fractional Sobolev embedding
$\dot W^{s, p}_\ast(\Omega)\hookrightarrow L^{\frac{np}{n-sp}}(\Omega)$), which means that 
for  any $u \in \dot{ W }_{*}^{s, p}(\Omega)$,
\begin{align}\label{ifsp}
\|u- u_{\Omega}\|_{L^{np/(n-sp)}(\Omega)} \leq C \|u\|_{\dot{ W }_{*}^{s, p}(\Omega)}  ,
\end{align}
holds, where $C$ depends on $n, s, p$ and $c$.
On the other hand,  additionally assume that $\Omega$ is a
bounded simply connected planar domain or a domain that is quasiconformally equaivalently to some uniform domain when $n\ge3$. They
\cite{d16, h13} also proved that   if    \eqref{ifsp} holds,
then $\Omega$ is a $c$-John domain.

If $1\leq p<\frac{n}{s}$,  it's easy to see that  $\phi_\frac{n}{s}(t)=Ct^\frac{np}{n-sp}$
 for any $t\ge0$,  where $C$ is a positive constant.
If  $ \phi(t)=t^p$ with $p\ge1 $ and $0<s<1$,
then the $(\phi_\frac{n}{s}, \phi)$-Poincar\'e inequality is the classical fractional $(\frac{np}{n-sp}, p)$-Poincar\'e inequality.

(iii) Analogue results to (ii) were established for  the intrinsic fractional Hajlasz-Sobolev space $\dot{M}_\ast^{s, p}(\Omega) $; see \cite{z11} for details.

\end{rem}

We also note that the imbeddings of the fractional Sobolev space
$\dot{ W}^{s, p}(\Omega)$  and
fractional Orlicz-Sobolev space $\dot{ V}^{s, \phi}(\Omega)$ were  were taken into account in the citations \cite{refc3, jw78,jw84,z123}.
Define  the fractional Orlicz-Sobolev space $\dot{ V }^{s, \phi}(\Omega)$   consisting of
  all functions $u\in L^1_{\loc}(\Omega) $   with
\begin{eqnarray*}
\|u\|_{\dot{ V }^{s, \phi}(\Omega)} := \inf \left\{\lambda > 0 : \int_{\Omega} \int_{\Omega}  \phi \left(\frac{|u(x)-u(y)|}{\lambda|x-y|^{s}}\right) \frac{dxdy}{|x-y|^{n}}\leq 1 \right\} < \infty.
\end{eqnarray*}
The $\dot{ V}^{s, \phi}(\Omega)$-(semi)norm is evidently derived by substituting  
the whole domain $\Omega$ for the range $B(x, \frac12\dist(x, \partial\Omega))$ for  the variable $y$
in  the $\dot{ V }^{s, \phi}_\ast(\Omega)$-(semi)norm.
 It goes without saying that $\dot{ V }^{s, \phi}(\rn) = \dot{ V}_{*}^{s, \phi}(\rn)$.
 For general domain $\Omega$, one always has $\dot{ V }^{s, \phi}(\Omega) \subset \dot{ V}_{*}^{s, \phi}(\Omega)$ with a normal bound,
 but the reverse side is not true necessarily.
When $\phi(t)=t^p$ with $p\ge1$, $\dot{ V}^{s, \phi}(\Omega)$ is  the
  fractional Sobolev space $\dot{W} ^{s, p}(\Omega) $, which  consists of
 all functions  $u\in L^1_\loc(\Omega)$ with
  \begin{eqnarray*}
\|u\|_{\dot{ W }^{s, p}(\Omega)} :=   \lf(\int_{\Omega} \int_\Omega \frac{|u(x)-u(y)|^p}{|x-y|^{n+sp}} \, dxdy\r)^{1/p}<\fz.
\end{eqnarray*}

\begin{rem}\rm
(i)  Let $s\in(0,1)$ and $1\le p<n/s$. It was shown in \cite{jw78,jw84,z123} that a domain  $\Omega$ supports the $\dot{ W }^{s, p}$-imbedding \begin{align*}
\|u-u_\Omega\|_{L^{\phi_\frac{n}{s}}(\Omega)} \leq C \|u\|_{\dot{V}^{s, \phi}(\Omega)} \quad \forall u \in  \dot{V}^{s, \phi}(\Omega).
\end{align*}
if and only if $\Omega$ is Ahlfors $n$-regular, that is,  there exists a constant $c>0$ such that
$$B(x,r)\cap\Omega|\ge Cr^n\quad\forall x\in\Omega,\ 0<r<2\diam\Omega.$$
Note that in he case $|\Omega|=\fz$   we set $u_\Omega=0$.

(ii)  Assume that $s\in(0,1)$ and Young function
 $\phi$ satisfies \eqref{eq2} and \eqref{eq3}.
It was shown in \cite{refc3} that
 Lipschitz domain $\Omega$  supports $\dot{ V}^{s, \phi}(\Omega)$-imbedding
\begin{align*}
\|u-u_\Omega\|_{L^{\phi_\frac{n}{s}}(\Omega)} \leq C \|u\|_{\dot{V}^{s, \phi}(\Omega)} \quad \forall u \in  \dot{V}^{s, \phi}(\Omega).
\end{align*}
But it is not clear  whether   Ahlfors $n$-regular domains characterize $\dot{ V}^{s, \phi}(\Omega)$-imbedding  domains.
\end{rem}

The paper is organized as follows.
The proof of  Theorem \ref{th1}(i) is given in section 2,  which  uses Boman's chain property,
 the embedding $\dot{ V }_\ast^{s, \phi}(Q) \hookrightarrow L^{\phi_{n/s}}(Q)$
for cubes $Q \subset \rn$
and the vector-valued inequality in Orlicz norms for the Hardy-Littlewood maximum operators.
We also give some property of $\phi \in \Delta_2$ with $K_\phi<2^{\frac{n}{s}}$  in section 2.
Conversely,  under the condition \eqref{eq1},  together with the aid of some ideas from \cite{bk95,  refh1,reflt3, refsh1, z123},  we
obtain the $LLC(2)$ property of $\Omega$,  and then
 prove Theorem \ref{th1}(ii) by a capacity argument;
see Section 3 for details.

\section{Proof of Theorem \ref{th1}(i) }

First we give the embedding $C_c^\infty(\Omega) \subset\dot{V}_\ast^{s, \phi}(\Omega)$.
It's easy to  know
\begin{equation}\label{eq1}
C_\phi :=\sup_{t>0}\int_0^t \frac{\phi(\rho)}{\phi(t)}\frac{d\rho}{\rho}<\infty.
\end{equation}
In fact, since for  practically all $t\ge0$ $\phi'(t)\ge0$ and $\phi'$ is increasing, we know
\begin{align*}
\frac{\phi(\rho)}{\rho}=\frac{\phi(\rho)-\phi(0)}{\rho}\leq \phi'(\rho).
\end{align*}
Hence
\begin{align*}
\int_0^t \frac{\phi(\rho)}{\phi(t)}\frac{d\rho}{\rho}
\leq \frac{1}{\phi(t)}\int_0^t\phi'(\rho) \,d\rho
\leq 1 ,
\end{align*}
that is, $C_\phi \leq 1$.

\begin{lem}\label{s2.l1}
Let $0 < s < 1$ , and $\phi$ be a  Young function satisfying $\eqref{eq1}$. For any bounded domain $\Omega \subset \mathbb{R}^n$,  we have $C_c^\infty(\Omega) \subset \dot{V}^{s, \phi}(\Omega)\subset\dot{V}_\ast^{s, \phi}(\Omega)$.
\end{lem}

\begin{proof}
$\forall u\in C_c^1(\Omega) $,  $L:=\|u\|_{L^\infty(\Omega)}+\|Du\|_{L^\infty(\Omega)}$, and $W\subset\Omega$ such that $V=\mathrm{supp } u\Subset W\Subset\Omega$,
then
 \begin{align*}
H:&=\int_\Omega\int_\Omega \phi\left(\frac{|u(x)-u(y)|}{\lambda |x-y|^s}\right)\frac{dxdy}{|x-y|^n}\\&
\leq\int_W\int_W \phi\left(\frac{L|x-y|}{\lambda |x-y|^s}\right)\frac{dxdy}{|x-y|^n}+2\int_V\int_{\Omega\setminus W} \phi\left(\frac{L}{\lambda |x-y|^s}\right)\frac{dxdy}{|x-y|^n}.
\end{align*}
By $\eqref{eq1}$, we have
\begin{align*}
\int_W\int_W \phi\left(\frac{L|x-y|}{\lambda |x-y|^s}\right)\frac{dxdy}{|x-y|^n}&\leq\int_W\int_{B(x, 2\diam W)}\phi\left(\frac{L|x-y|^{1-s}}{\lambda }\right)\frac{dy}{|x-y|^n}dx\\&
=n\omega_n\int_W\int_0^{2\diam W}\phi\left(\frac{L\rho^{1-s}}{\lambda }\right)\frac{d\rho}{\rho}dx\\&
=n\omega_n\frac{1}{1-s}\int_W\int_0^{\frac{L(2\diam W)^{1-s}}{\lambda}}\phi(\mu)\frac{d\mu}{\mu}dx\\&
\leq C_\phi n\omega_n\frac{1}{1-s}\phi\left(\frac{L(2\diam W)^{1-s}}{\lambda}\right)|W|.
\end{align*}
And
\begin{align*}
\int_V\int_{\Omega\setminus W} \phi\left(\frac{L}{\lambda |x-y|^s}\right)\frac{dxdy}{|x-y|^n}&
\leq\int_V\int_{\Omega\setminus B(y, \dist(V, W^{\complement}))}\phi\left(\frac{L}{\lambda |x-y|^s}\right)\frac{dx}{|x-y|^n}dy\\&
\leq n\omega_n \int_V\int_{\dist(V, W^{\complement})}^\infty \phi\left(\frac{L}{\lambda \rho^s}\right)\frac{d\rho}{\rho}dy\\&
=n\omega_n\frac{1}{s}\int_V\int_0^{\frac{L}{\lambda \dist(V, W^{\complement})^s}}\phi(\mu)\frac{d\mu}{\mu}dy\\&
\leq C_\phi n\omega_n\frac{1}{s}\phi\left(\frac{L}{\lambda \dist(V, W^{\complement})^s}\right)|V|.
\end{align*}
If $\lambda$ is so large,  we have  $H\leq1$,with $u\in\dot{V}^{s,\phi}(\Omega)$,
so $C_c^1(\Omega) \subset \dot{V}^{s,\phi}(\Omega)$.
 Combining $C_c^\infty(\Omega) \subset  C_c^1(\Omega)$ and $\dot{V}^{s,\phi}(\Omega)\subset\dot{V}_\ast^{s,\phi}(\Omega)$, we get the result desired.
\end{proof}

To prove Theorem \ref{th1}(i),
we need  the embedding $\dot{V}^{s,\phi}(Q)\hookrightarrow L^{\phi_{n/s}}(Q)$ in all cubes $Q\subset\rn$.
So firstly, we give some lemmas we needed.
\begin{lem}\label{q1}
 Let $\phi \in \Delta_2$ be a Young function, then $\forall c>1,\phi(cx)\leq c^{K_\phi-1}\phi(x)$.
\end{lem}
\begin{proof}
By the increasing property of $\phi^\prime$,
$$\phi(2x)-\phi(x)=\int_x^{2x}\phi^\prime(t)dt\geq \phi^\prime(x)  ,~\forall x>0.$$
$\phi \in \Delta_2, \phi(2x)-\phi(x)\leq(K_\phi-1)\phi(x)$,
so
$$(\ln \phi)^\prime(x)=\frac{\phi^\prime(x)}{\phi(x)}\leq\frac{K_\phi-1}{x}.$$
For any $c>1$,we have
$$\ln\left(\frac{\phi(cx)}{\phi(x)}\right)=\int_x^{cx}(\ln\phi)^\prime(t)dt\leq\int_x^{cx}\frac{K_\phi-1}{t}dt
=\ln(c^{K_\phi-1}).$$
So $\phi(cx)\leq c^{K_\phi-1}\phi(x)$.
\end{proof}
\begin{lem}\label{le1}
 Let $\phi \in \Delta_2$ be a Young function satisfying  $K_\phi<2^{\frac{n}{s}}$, then $\phi$ satisfies $\eqref{eq2}$ and $\eqref{eq3}$.
\end{lem}
\begin{proof}
By the definition of the $K_\phi$ in \eqref{da2}, we get $\phi(2t)\leq K_\phi\phi(t)$. Hence
\begin{align*}
\int_\frac{t}{2}^t\left(\frac{\tau}{\phi(\tau)}\right)^\frac{s}{n-s}d\tau&=\int_\frac{t}{4}^\frac{t}{2}
\left(\frac{2\tau}{\phi(2\tau)}\right)^\frac{s}{n-s}2d\tau
\\&\geq\int_\frac{t}{4}^\frac{t}{2}\left(\frac{2\tau}{K_\phi\phi(\tau)}\right)^\frac{s}{n-s}2d\tau.
\end{align*}
Then
\begin{align*}
\int_\frac{t}{4}^\frac{t}{2}\left(\frac{\tau}{\phi(\tau)}\right)^\frac{s}{n-s}d\tau
\leq\frac{K_\phi^\frac{s}{n-s}}{2^\frac{n}{n-s}}\int_\frac{t}{2}^t\left(\frac{\tau}{\phi(\tau)}\right)
^\frac{s}{n-s}d\tau.
\end{align*}

By induction, we have
\begin{align*}
\int_\frac{t}{2^m}^\frac{t}{2^{m-1}}\left(\frac{\tau}{\phi(\tau)}\right)^\frac{s}{n-s}d\tau
&\leq\frac{K_\phi^\frac{s}{n-s}}{2^\frac{n}{n-s}}\int_\frac{t}{2^{m-1}}^\frac{t}{2^{m-2}}\left(\frac{\tau}
{\phi(\tau)}\right)^\frac{s}{n-s}d\tau
\\&\leq\left(\frac{K_\phi^\frac{s}{n-s}}{2^\frac{n}{n-s}}\right)^{m-1}
\int_\frac{t}{2}^t\left(\frac{\tau}{\phi(\tau)}\right)^\frac{s}{n-s}d\tau.
\end{align*}
Change $m$ from 1 to $\infty$ and sum up, we can get
\begin{equation*}
\int_0^t\left(\frac{\tau}{\phi(\tau)}\right)^\frac{s}{n-s}d\tau\leq\sum_{m=1}^\infty \left(\frac{K_\phi^\frac{s}{n-s}}{2^\frac{n}{n-s}}\right)^{m-1}\int_\frac{t}{2}^t\left(\frac{\tau}
{\phi(\tau)}\right)^\frac{s}{n-s}d\tau.
\end{equation*}
The series convergences because if the range of the $K_\phi$,and
\begin{align*}
\left(\frac{t}{\phi(t)}\right)^\prime(t)&=\frac{\phi(t)-t\phi^\prime(t)}{\phi^2(t)}
\\&=\frac{\frac{\phi(t)-\phi(0)}{t}-\phi^\prime(t)}{\phi^2(t)}\\&=
\frac{\phi^\prime(\xi)-\phi^\prime(t)}{\phi^2(t)},
\end{align*}
where $0<\xi<t$, by the convexity of  $\phi$, we know $\left(\frac{t}{\phi(t)}\right)^\prime(t)\leq0$, then \eqref{eq2} follows from decreasing property of $\frac{\tau}{\phi(\tau)}$,
actually,
\begin{equation*}
\int_\frac{t}{2}^t\left(\frac{\tau}{\phi(\tau)}\right)^\frac{s}{n-s}d\tau\leq\left(\frac{\frac{t}{2}}
{\phi(\frac{t}{2})}\right)^\frac{s}{n-s}\frac{t}{2}<\infty.
\end{equation*}
Similarly,
\begin{align*}
\int_0^{2^m}\left(\frac{\tau}{\phi(\tau)}\right)^\frac{s}{n-s}d\tau&\geq\int_0^{2^{m-1}}\left(\frac{2\tau}
{K_\phi\phi(\tau)}\right)^\frac{s}{n-s}2d\tau\\&=\frac{2^\frac{n}{n-s}}{K_\phi^\frac{s}{n-s}}\int_0^{2^{m-1}}
\left(\frac{\tau}{\phi(\tau)}\right)^\frac{s}{n-s}d\tau\\&\geq\ldots\\&\geq\left(\frac{2^\frac{n}{n-s}}
{K_\phi^\frac{s}{n-s}}\right)^{m}\int_0^1\left(\frac{\tau}{\phi(\tau)}\right)^\frac{s}{n-s}d\tau, ~\forall m \in \mathbb{N}.
\end{align*}
Let $m \to \infty$ we get \eqref{eq3}.
\end{proof}

\begin{lem}\label{te}
 Let $\phi \in \Delta_2$ be a Young function satisfying  $K_\phi<2^{\frac{n}{s}}$, then
\begin{equation}\label{te1}
\frac{H(A)}{A}\leq C\frac{1}{{\phi(A)}^\frac{s}{n}}.
\end{equation}
\end{lem}
\begin{proof}
By Lemma \ref{le1}, we have
\begin{equation*}
\int_0^t\left(\frac{\tau}{\phi(\tau)}\right)^\frac{s}{n-s}d\tau\leq C\int_\frac{t}{2}^t\left(\frac{\tau}{\phi(\tau)}\right)^\frac{s}{n-s}d\tau
\leq C(\frac{\frac{t}{2}}{\phi(\frac{t}{2})})^\frac{s}{n-s}\frac{t}{2}.
\end{equation*}
Then
\begin{align*}
\frac{H(A)}{A}&=\frac{\left(\int_0^A(\frac{\tau}{\phi(\tau)})^\frac{s}{n-s}d\tau\right)^\frac{n-s}{n}}{A}
\leq \frac{\left(C(\frac{\frac{A}{2}}{\phi(\frac{A}{2})})^\frac{s}{n-s}\frac{A}{2}\right)^\frac{n-s}{n}}{A}
\leq \frac{\left(C\left(\frac{\frac{A}{2}}{\frac{1}{K_\phi}\phi(A)}\right)^\frac{s}{n-s}\frac{A}{2}\right)
^\frac{n-s}{n}}{A}
\leq C\frac{1}{{\phi(A)}^\frac{s}{n}}.
\end{align*}
\end{proof}

With above lemmas,  we proved $\dot{V}^{s, \phi}(Q)\hookrightarrow L^{\phi_{n/s}}(Q)$.
\begin{lem}\label{le2}
Let $0 < s < 1$ and $\phi \in \Delta_2$ be a Young function satisfying $K_\phi<2^{\frac{n}{s}}$,  then there exists a constant $C_1 = C_1(n, s)$ such that

\begin{equation}\label{eq2.1}
\int_Q \phi_\frac{n}{s}\left(\frac{u(x)-u_Q}{\lambda}\right)dx\leq\int_Q\int_Q \phi\left(\frac{C|u(x)-u(y)|}{\lambda{|x-y|}^s}\right)\frac{dxdy}{{|x-y|}^n}.
\end{equation}
for all cubes $Q \subset \mathbb{R}^n,  ~ u \in {\dot {V}}^{s, \phi}(Q)$ and $\lambda \geq C_1\|u\|_{{\dot {V}}^{s, \phi}(Q)}$.
\end{lem}
\begin{proof}
Denote a cube centered at the origin with sides of length 2 paralleled to the axes by $Q(0,1)$.
At first we prove that
\begin{equation}
\int_{Q(0, 1)}  \phi_\frac{n}{s}\left(\frac{|u(x)-u_{Q(0, 1)} |}{\lambda}\right)dx\leq\int_{Q(0, 1)} \int_{Q(0, 1)}  \phi\left(\frac{C_1|u(x)-u(y)|}{\lambda{|x-y|}^s}\right)\frac{dxdy}{{|x-y|}^n},
\end{equation}
where $u\in {\dot {V}}^{s, \phi}({Q(0, 1)}), \, \lambda \geq C_1\|u\|_{{\dot {V}}^{s, \phi}(Q(0, 1))}$.

$K_\phi<2^{\frac{n}{s}}$, by Lemma \ref{le1} and \cite{refc3}, we have
\begin{equation*}
\|{u}\|_{L^{\phi_\frac{n}{s}}({Q(0, 1)})}\leq {C_1}\|{u}\|_{{\dot {V}}^{s, \phi}({Q(0, 1)})}.
\end{equation*}
where
$$u\in {\dot {V}}^{s, \phi}_{\perp}({Q(0, 1)}), C_1=C_1(n, s), $$
$${\dot {V}}^{s, \phi}_{\perp}({Q(0, 1)}):=\left\{ u\in{\dot {V}}^{s, \phi}({Q(0, 1)}):u_{Q(0, 1)}=0\right\}.$$
Replacing $u$ by $u-u_{Q(0, 1)}$,  we have
\begin{equation*}
\|{u-u_{Q(0, 1)}}\|_{L^{\phi_\frac{n}{s}}({Q(0, 1)})}\leq {C_1}\|{u-u_{Q(0, 1)}}\|_{{\dot {V}}^{s, \phi}({Q(0, 1)})},
\end{equation*}
where $u\in {\dot {V}}^{s, \phi}({Q(0, 1)})$.
When $\|u\|_{{\dot {V}}^{s, \phi}(Q(0, 1))}=0$, $u$ is constant in $Q(0, 1)$, the equality holds. Suppose that $\|u\|_{{\dot {V}}^{s, \phi}(Q(0, 1))}\neq0$,
then
$$\int_{Q(0, 1)}\phi_\frac{n}{s}\left(\frac{|u-u_{Q(0, 1)}|}{C_1\|u\|_{{\dot {V}}^{s, \phi}({Q(0, 1)})}}\right)dx\leq1.$$
Fixed $u_0 \in {\dot {V}}^{s, \phi}({Q(0, 1)})$,  let
$$M:=\int_{Q(0, 1)}\int_{Q(0, 1)}\phi\left(\frac{C_1|u_0(x)-u_0(y)|}{|x-y|^s}\right)\frac{dxdy}{|x-y|^n}\neq0.$$
Let $\overline{\phi}=\frac{\phi}{M}$, then $\overline{\phi}_\frac{n}{s}(t)=\frac{1}{M}\phi_\frac{n}{s}(\frac{t}{M^\frac{s}{n}})$ and $C_1=C_1(n, s)$,
so $$\|u-u_{Q(0, 1)}\|_{L^{\overline{\phi}_\frac{n}{s}}({Q(0, 1)})}\leq C_1\|u\|_{{\dot {V}}^{s, \overline{\phi}}({Q(0, 1)})}.$$
Then we get $$\int_{Q(0, 1)}\phi_\frac{n}{s}\left(\frac{|u-u_{Q(0, 1)}|}{C_1M^\frac{s}{n}\|u\|_{{\dot {V}}^{s, \overline{\phi}}({Q(0, 1)})}}\right)dx\leq M.$$
And $C_1\|u_0\|_{{\dot {V}}^{s, \overline{\phi}}({Q(0, 1)})}\leq 1$, otherwise,
\begin{align*}
1&<\int_{Q(0, 1)}\int_{Q(0, 1)}\overline{\phi}\left(\frac{C_1|u_0(x)-u_0(y)|}{|x-y|^s}\right)
\frac{dxdy}{|x-y|^n}\\
&=\frac{1}{M}\int_{Q(0, 1)}\int _{Q(0, 1)}{\phi}\left(\frac{C_1|u_0(x)-u_0(y)|}{|x-y|^s}\right)\frac{dxdy}{|x-y|^n}=1,
\end{align*}
we get a contradiction.

Specially,  when $u=u_0$, we have
\begin{align*}
&\int_{Q(0, 1)}\phi_\frac{n}{s}\left(\frac{|u_0-u_{0{Q(0, 1)}}|}{(\int_{Q(0, 1)}\int_{Q(0, 1)}
\phi\left(\frac{C_1|u_0(x)-u_0(y)|}{|x-y|^s}\right)\frac{dxdy}{|x-y|^n})^\frac{s}{n}}\right)dx\\
\leq&\int_{Q(0, 1)}\phi_\frac{n}{s}\left(\frac{|u_0-u_{0{Q(0, 1)}}|}{C_1M^\frac{s}{n}\|u_0\|_{{\dot {V}}^{s, \overline{\phi}}({Q(0, 1)})}}\right)dx
\\\leq & \int_{Q(0, 1)}\int_{Q(0, 1)}\phi\left(\frac{C_1|u_0(x)-u_0(y)|}{|x-y|^s}\right)\frac{dxdy}{|x-y|^n}.
\end{align*}
By the arbitrariness of $u_0$, we have
\begin{align*}
&\int_{Q(0, 1)}\phi_\frac{n}{s}\left(\frac{|u-u_{{Q(0, 1)}}|}{(\int_{Q(0, 1)}\int_{Q(0, 1)}
\phi\left(\frac{C_1|u(x)-u(y)|}{|x-y|^s}\right)\frac{dxdy}{|x-y|^n})^\frac{s}{n}}\right)dx\\
\leq & \int_{Q(0, 1)}\int_{Q(0, 1)}\phi\left(\frac{C_1|u(x)-u(y)|}{|x-y|^s}\right)\frac{dxdy}{|x-y|^n}.
\end{align*}
Replacing $u$ by $\frac{u}{\lambda}$,
\begin{align}\label{sadfs}
&\int_{Q(0, 1)}\phi_\frac{n}{s}\left(\frac{|u-u_{{Q(0, 1)}}|}{\lambda(\int_{Q(0, 1)}\int_{Q(0, 1)}
\phi\left(\frac{C_1|u(x)-u(y)|}{\lambda|x-y|^s}\right)\frac{dxdy}{|x-y|^n})^\frac{s}{n}}\right)dx\\
\leq & \int_{Q(0, 1)}\int_{Q(0, 1)}\phi\left(\frac{C_1|u(x)-u(y)|}{\lambda|x-y|^s}\right)\frac{dxdy}{|x-y|^n}.
\end{align}

Let $\lambda\geq C_1\|u\|_{{\dot {V}}^{s, \phi}({Q(0, 1)})}$, then$$\int_{Q(0, 1)}\int_{Q(0, 1)}\phi\left(\frac{C_1|u(x)-u(y)|}{\lambda|x-y|^s}\right)\frac{dxdy}{|x-y|^n}\leq1, $$

so
\begin{align*}
\int_{Q(0, 1)}\phi_\frac{n}{s}\left(\frac{|u-u_{{Q(0, 1)}}|}{\lambda}\right)dx&\leq\int_{Q(0, 1)}\phi_\frac{n}{s}\left(\frac{|u-u_{{Q(0, 1)}}|}{\lambda(\int_{Q(0, 1)}\int_{Q(0, 1)}\phi\left(\frac{C_1|u(x)-u(y)|}{\lambda|x-y|^s}\right)\frac{dxdy}{|x-y|^n})^\frac{s}{n}}\right)dx
\\&\leq\int_{Q(0, 1)}\int_{Q(0, 1)}\phi\left(\frac{C_1|u(x)-u(y)|}{\lambda|x-y|^s}\right)\frac{dxdy}{|x-y|^n}.
\end{align*}

Now we prove the case of general cube $Q$.
Let $Q$ be a cube with $a$ as the center and $2l$ as the side length,  then there is an orthogonal transformation $T$, and $T(Q-a)=Q(0, l)$.
$\forall u \in {\dot {V}}^{s, \phi}(Q)$ and $u$ is not a constant.
Let $v(x)=\frac{u(T^{-1}(lx)+a)}{ l^s}, $ where $x\in Q(0, 1)$, then $v\in {\dot {V}}^{s, \phi}(Q(0, 1))$
 and $v_{Q(0, 1)}=\frac{u_Q}{l^s}$.
And
\begin{align*}
&\int_{Q(0, 1)}\int_{Q(0, 1)}\phi\left(\frac{C_1|v(x)-v(y)|}{\lambda|x-y|^s}\right)\frac{dxdy}{|x-y|^n}\\
&=\int_{Q(0, 1)}\int_{Q(0, 1)}\phi\left(\frac{C_1|\frac{u(T^{-1}(lx)+a)}{ l^s}-\frac{u(T^{-1}(ly)+a)}{ l^s}|}{\lambda|x-y|^s}\right)\frac{dxdy}{|x-y|^n},
\end{align*}
by transformation $z_1=T^{-1}(lx)+a, z_2=T^{-1}(ly)+a, $
we have $|x-y|=\left|\frac{T(z_1-a)}{l}-\frac{T(z_2-a)}{l}\right|=\frac{|z_1-z_2|}{l}, $
so\begin{align*}
\int_{Q(0, 1)}\int_{Q(0, 1)}\phi\left(\frac{C_1|v(x)-v(y)|}{\lambda|x-y|^s}\right)\frac{dxdy}{|x-y|^n}
=\int_Q\int_Q\phi\left(\frac{C_1|{u(z_1)}-{u(z_2)}|}{\lambda|z_1-z_2|^s}\right)\frac{dz_1dz_2}{l^n|z_1-z_2|^n},
\end{align*}
and
\begin{align*}
&\int_{Q(0, 1)}\phi_\frac{n}{s}\left(\frac{|v-v_{{Q(0, 1)}}|}{\lambda(\int_{Q(0, 1)}\int_{Q(0, 1)}
\phi\left(\frac{C_1|u(x)-u(y)|}{\lambda|x-y|^s}\right)\frac{dxdy}{|x-y|^n})^\frac{s}{n}}\right)dx\\
&=\int_{Q(0, 1)}\phi_\frac{n}{s}\left(\frac{|v-v_{{Q(0, 1)}}|}{\lambda
(\int_Q\int_Q\phi\left(\frac{C_1|{u(z_1)}-{u(z_2)}|}{\lambda|z_1-z_2|^s}\right)
\frac{dz_1dz_2}{l^n|z_1-z_2|^n})^\frac{s}{n}}\right)dx.
\end{align*}
By transformation $y=T^{-1}(lx)+a$, we get
\begin{align*}
&\int_{Q(0, 1)}\phi_\frac{n}{s}\left(\frac{|v-v_{{Q(0, 1)}}|}{\lambda(\int_{Q(0, 1)}\int_{Q(0, 1)}\phi
\left(\frac{C_1|u(x)-u(y)|}{\lambda|x-y|^s}\right)\frac{dxdy}{|x-y|^n})^\frac{s}{n}}\right)dx\\
&=\int_{Q}\phi_\frac{n}{s}\left(\frac{|u(y)-u_{{Q}}|}{\lambda(\int_Q\int_Q\phi
\left(\frac{C_1|{u(z_1)}-{u(z_2)}|}{\lambda|z_1-z_2|^s}\right)
\frac{dz_1dz_2}{|z_1-z_2|^n})^\frac{s}{n}}\right)\frac{dy}{l^n}.
\end{align*}
By \eqref{sadfs}, we have
\begin{align*}
&\int_{Q}\phi_\frac{n}{s}\left(\frac{|u(y)-u_{{Q}}|}
{\lambda(\int_Q\int_Q\phi\left(\frac{C_1|{u(z_1)}-{u(z_2)}|}{\lambda|z_1-z_2|^s}\right)
\frac{dz_1dz_2}{|z_1-z_2|^n})^\frac{s}{n}}\right){dy}\\
\leq& \int_Q\int_Q\phi\left(\frac{C_1|{u(z_1)}-{u(z_2)}|}{\lambda|z_1-z_2|^s}\right)
\frac{dz_1dz_2}{|z_1-z_2|^n}.
\end{align*}
Let $\lambda\geq C_1\|u\|_{{\dot {V}}^{s, \phi}({Q})}$, then$$\int_{Q}\int_{Q}\phi\left(\frac{C_1|u(x)-u(y)|}{\lambda|x-y|^s}\right)\frac{dxdy}{|x-y|^n}\leq1.$$
Hence,
\begin{align*}
\int_{Q}\phi_\frac{n}{s}\left(\frac{|u-u_{Q}|}{\lambda}\right)dx&\leq\int_{Q}\phi_\frac{n}{s}
\left(\frac{|u-u_{Q}|}{\lambda(\int_Q\int_Q\phi\left(\frac{C_1|u(x)-u(y)|}{\lambda|x-y|^s}\right)
\frac{dxdy}{|x-y|^n})^\frac{s}{n}}\right)dx
\\&\leq\int_Q\int_Q\phi\left(\frac{C_1|u(x)-u(y)|}{\lambda|x-y|^s}\right)\frac{dxdy}{|x-y|^n}.
\end{align*}
\end{proof}

We also  need the Fefferman-Stein type verct-valued inequality for Hardy-Littlewood maximum operator in Orlicz space. Denote by $\mathcal{M}$ the Hardy-Littlewood maximum operator,
\begin{equation*}
\mathcal{M}(g)(x)=\sup_{x\in Q}\fint_{Q}|g|dx
\end{equation*}
with the supremum taken over all cubes $Q \subset \mathbb{R}^n$ containing x.
The Young function $\phi$ is in  $\nabla_{2}$ if there exist a $a>1$,  such that
\begin{equation*}
\phi(x)\leq\frac{1}{2a}\phi(ax), ~\forall x\geq0.
\end{equation*}
\begin{lem}\label{le2.3}
If $\phi \in \Delta_2$ be a Young function satisfying $K_{\phi}<2^{\frac{n}{s}}$, then $\phi_\frac{n}{s} \in \Delta_2\cap\nabla_2$.
\end{lem}
\begin{proof}
We know
\begin{align*}
H(2t)&=\left(\int_0^{2t}(\frac{\tau}{\phi(\tau)})^\frac{s}{n-s}d\tau\right)^\frac{n-s}{n}\\
&\geq\left(\int_0^t(\frac{2\tau}{K_{\phi}\phi(\tau)})^\frac{s}{n-s}2d\tau\right)^\frac{n-s}{n}
=\frac{2}{K_{\phi}^\frac{s}{n}}H(t).
\end{align*}

Letting $2y=H(2t)$,  we have $K_{\phi}^\frac{s}{n}y\geq H\left(\frac{H^{-1}(2y)}{2}\right)$.
Therefore,
$$H^{-1}(2y)\leq2H^{-1}(K_{\phi}^\frac{s}{n}y)\leq2^2H^{-1}(K_{\phi}^\frac{s}{n}\frac {K_{\phi}^\frac{s}{n}}{2}y)\leq...\leq2^{m+1}H^{-1}(K_{\phi}^\frac{s}{n}\left(\frac {K_{\phi}^\frac{s}{n}}{2}\right)^my).$$
Because of the range of K,  we get $\frac{K_{\phi}^\frac{s}{n}}{2}<1$.
Let $m$ so big that $K_{\phi}^\frac{s}{n}\left(\frac {K_{\phi}^\frac{s}{n}}{2}\right)^m<1$.
Then we have $H^{-1}(2y)<CH^{-1}(y)$.
So  $H^{-1}\in \Delta_2$  and $\phi_\frac{n}{s}=\phi \circ H^{-1}\in \Delta_2$.

By the decreasing property of $\frac{\tau}{\phi(\tau)}$,
\begin{align*}
H(2^\frac{n}{s}x)&=\left(\int_0^{2^\frac{n}{s}x}\left(\frac{\tau}{\phi(\tau)}\right)^\frac{s}{n-s}d\tau\right)
^\frac{n-s}{n}
\\&=\left(\int_0^x\left(\frac{2^\frac{n}{s}\tau}{\phi(2^\frac{n}{s}\tau)}\right)^\frac{s}{n-s}
2^\frac{n}{s}d\tau
\right)^\frac{n-s}{n}\\
&\leq\left(\int_0^x\left(\frac{\tau}{\phi(\tau)}\right)^\frac{s}{n-s}2^\frac{n}{s}d\tau\right)^\frac{n-s}{n}
\\&=2^\frac{n-s}{s}H(x).
\end{align*}
So ~$2^\frac{n}{s}x\leq H^{-1}(2^\frac{n-s}{s}H(x))$, ~ then~ $2^\frac{n}{s}H^{-1}(x)\leq H^{-1}(2^\frac{n-s}{s}x)$.

And we have
 \begin{align*}
2^\frac{n}{s}\phi \circ H^{-1}(x)\leq \phi(2^\frac{n}{s}H^{-1}(x))\leq\phi \circ H^{-1}(2^\frac{n-s}{s}x).
\end{align*}
Letting $a=2^\frac{n-s}{s}>1$, we have $\phi_\frac{n}{s}(x)\leq\frac{1}{2a}\phi_\frac{n}{s}(ax)$ and $\phi_\frac{n}{s}\in \nabla_2$.
\end{proof}

\begin{rem}
If $K_\phi\geq2^{\frac{n}{s}}$, there exists $\phi\in\Delta_2$ such that $\phi\frac{n}{s}\notin\Delta_2$.\cite{refc3} Example 6.4:
Let $\phi$ with
\begin{align*}
\phi(t)=\left\{\begin{array}{ll}
{t^\frac{n}{s}(\log\frac{1}{t})^{\alpha_0}}& near ~ zero , \\
{t^ \frac{n}{s}(\log{t})^{\alpha}}& near ~infinity,
\end{array}\right.
\end{align*}
where $\alpha_0>\frac{n}{s}-1, \alpha\leq\frac{n}{s}-1$.And connected by a convex function,  then
\begin{align*}
\phi_\frac{n}{s}(t)~is ~equivalent ~to \left\{\begin{array}{ll}
{e^{-t^{-\frac{n}{s(\alpha_0+1)-n}}}}& near ~zero, \\
{e^{t^{\frac{n}{n-s(\alpha+1)}}}}& near ~infinity, \alpha<\frac{n}{s}-1, \\{e^{e^{t^\frac{n}{n-s}}}}&near~ infinity, \alpha=\frac{n}{s}-1,
\end{array}\right.
\end{align*}
so $\phi_\frac{n}{s}\notin\Delta_2$.
\end{rem}

We then propose a few lemmas that might be utilized to support the assertion of  Theorem \ref{th1}(i).

\begin{lem}[\cite{cmp11}]\label{le2.4}
Let $\psi \in \Delta_2\cap\nabla_2$ be a Young function. For any $0<q<\infty$,  there exists a constant $C>1$ depending on $n,  q,  K_\psi$ and $a$ such that for all sequences $\{f_j\}_{j\in\mathbb{N}}$ , we have
\begin{equation*}
\int_{\mathbb{R}^n}\psi\left(\left[\sum_{j\in\mathbb{N}}(\mathcal{M}(f_j))^2\right]^\frac{1}{q}\right) dx\leq C(n, K_\psi, a)\int_{\mathbb{R}^n}\psi\left(\left[\sum_{j\in\mathbb{N}}(f_j)^2\right]^\frac{1}{q}\right) dx.
\end{equation*}
\end{lem}
\begin{lem}\label{le2.5}
For any constant $k\geq1$, sequence $\{a_j\}_{j\in\mathbb{N}}$, and cubes $\{Q_j\}_{j\in\mathbb{N}}$ with $\sum_{j}\chi_{Q_j}\leq k$, we have
\begin{equation*}
\sum_{j}|a_j|\chi_{kQ_j}\leq C(k, n)\sum_j[\mathcal{M}(|a_j|^{\frac{1}{2}}\chi_{Q_j})]^2.
\end{equation*}
\end{lem}
\begin{proof}
By the definition of $\mathcal{M}$, ~we know
\begin{equation*}
\chi_{kQ_j}\leq k^n\mathcal{M}({\chi}_{Q_j}).
\end{equation*}
So
\begin{equation*}
\sum_{j}|a_j|\chi_{kQ_j}=\sum_{j}(|a_j|^\frac{1}{2}\chi_{kQ_j})^2\leq k^{2n}\sum_j[\mathcal{M}(|a_j|^{\frac{1}{2}}\chi_{Q_j})]^2.
\end{equation*}
\end{proof}

Now let us begin to give the proof of Theorem \ref{th1}(i).
\begin{proof}[Proof of Theorem \ref{th1}(i)]
Let $\Omega$ be a c-John domain. By Boman   \cite{b82} and Buckley \cite{bkl95},  $\Omega$ enjoys the following chain property: for every integer $\kappa>1$,  there exist a positive constant $C(\kappa,  \Omega)$ and a collection $\mathcal{F}$ of the cubes such that

(i) $Q\subset \kappa Q \subset \Omega$  for all $Q \in \mathcal{F},  \Omega =\cup_{Q\in \mathcal{F}}Q$ and
$$\sum _{Q\in\mathcal{F}}\chi_{\kappa Q}  \leq  C_{\kappa, c} \chi_{\Omega}. $$

(ii) $Q_{0} \in \mathcal{ F} $ is a fixed cube. For any other $Q \in  \mathcal{F}$,
there exist a subsequence $\{Q_{j}\}_{j=1}^{N} \subset \mathcal{F}$,  satisfying that $Q = Q_{N} \subset C_{\kappa, c}  Q_{j} $,  $C_{\kappa, c}^{-1}|Q_{j+1}| \leq |Q_{j} | \leq C_{\kappa, c}|Q_{j+1}|$ and
 $|Q_{j} \cap Q_{j+1}|\geq C_{\kappa, c}^{-1} \min\{|Q_{j} |,  |Q_{j+1}|\}$ for all $j =0,  \ldots,  N-1$.

Let $\kappa=5n$,  by (i) $Q\subset 5nQ\subset \Omega$ for each $Q\in \mathcal{F}$,
$$d(Q, \partial\Omega)\geq d(Q, \partial(5nQ)\geq\frac{5n-1}{2}l(Q)\geq2nl(Q), $$and hence
$$|x-y|\leq\sqrt{n}l(Q)\leq nl(Q)\leq \frac{1}{2}d(Q, \partial\Omega)\leq\frac{1}{2}d(x, \partial\Omega), ~\forall x, y\in Q\in\mathcal{F}.$$
Let $u\in \dot{V}_{*}^{s, \phi}(\Omega)$.
Up to  approximating by $\min\{\max\{u, -N\}, N\}$,  we can assume that $u\in L^\infty(\Omega)$,  and by the boundedness of $\Omega$,   $u\in L^1(\Omega)$.

By $$\forall x, y \in Q,  |x-y|\leq\frac{1}{2}d(x, \partial\Omega), $$we know
$$\int_Q\int_Q\phi\left(\frac{|u(x)-u(y)|}{{{\lambda}}{|x-y|}^s}\right)\frac{dxdy}{{|x-y|}^n}
\leq\int_Q\int_{B(x, \frac{1}{2}d(x, \partial\Omega))}\phi\left(\frac{|u(x)-u(y)|}{{{\lambda}}{|x-y|}^s}\right)\frac{dydx}{{|x-y|}^n}, $$
then $\|u\|_{{\dot {V}}^{s, \phi}(Q)}\leq\|u\|_{{\dot {V}_\ast}^{s, \phi}(\Omega)}, $
so
\begin{equation}\label{gjs}
\lambda\geq\|u\|_{{\dot {V}}^{s, \phi}(Q)}~ when~ \lambda\geq\|u\|_{{\dot {V}_\ast}^{s, \phi}(\Omega)}.
\end{equation}
Because of the convexity of $\phi_\frac{n}{s}$,  we have
\begin{align*}
I:&=\int_{\Omega} \phi_{\frac{n}{s}}\left(\frac{|u(z)- u_{\Omega}|}{\lambda}\right)dz\\&
\leq\int_{\Omega} \phi_{\frac{n}{s}}\left(\frac{1}{2}\left(\frac{2|u(z)- u_{Q_0}|+2|u_{\Omega}- u_{Q_0}|}{\lambda}\right)\right)dz\\&
\leq \frac{1}{2}\left[\int_{\Omega} \phi_{\frac{n}{s}}\left(\frac{2|u(z)- u_{Q_0}|}{\lambda}\right)dz+|\Omega|\phi_{\frac{n}{s}} \left(\frac{2|u_{\Omega}- u_{Q_0}|}{\lambda}\right)\right].
\end{align*}
By Jensen inequality, $$|\Omega|\phi_{\frac{n}{s}} \left(\frac{2|u_{\Omega}- u_{Q_0}|}{\lambda}\right)\leq\int_{\Omega} \phi_{\frac{n}{s}}\left(\frac{2|u(z)- u_{Q_0}|}{\lambda}\right)dz.$$
In (i) we have $\chi_\Omega\leq\sum_{Q\in \mathcal{F}}\chi_Q$,  so
\begin{align*}
I&\leq\int_{\Omega} \phi_{\frac{n}{s}}\left(\frac{2|u(z)- u_{Q_0}|}{\lambda}\right)dz\\&
\leq\sum_{Q\in\mathcal{F}}\int_Q\phi_{\frac{n}{s}}\left(\frac{2|u(z)- u_{Q_0}|}{\lambda}\right)dz\\&
\leq\frac{1}{2}\sum_{Q\in\mathcal{F}}\int_Q\phi_{\frac{n}{s}}\left(\frac{4|u(z)- u_{Q}|}{\lambda}\right)dz+\frac{1}{2}\sum_{Q\in\mathcal{F}
\setminus\{Q_0\}}|Q|\phi_{\frac{n}{s}} \left(\frac{4|u_{Q}- u_{Q_0}|}{\lambda}\right)\\&:=\frac{1}{2}I_1+\frac{1}{2}I_2.
\end{align*}
By the inequality \eqref{eq2.1}, \eqref{gjs} and $$\forall x, y \in Q,  |x-y|\leq\frac{1}{2}d(x, \partial\Omega), $$ we know
\begin{align*}
I_1&\leq\sum_{Q\in\mathcal{F}}\int_Q\int_Q\phi\left(\frac{|u(x)-u(y)|}{{\frac{\lambda}{4C_1}}{|x-y|}^s}\right)
\frac{dxdy}{{|x-y|}^n}\\&
\leq\sum_{Q\in\mathcal{F}}\int_Q\int_{B(x, \frac{1}{2}d(x, \partial\Omega))}\phi\left(\frac{|u(x)-u(y)|}{{\frac{\lambda}{4C_1}}{|x-y|}^s}\right)\frac{dydx}{{|x-y|}^n}.
\end{align*}
Using the $\sum _{Q\in\mathcal{F}}\chi_{\kappa Q}  \leq  C_{\kappa, c} \chi_{\Omega}$ in (i) above,
\begin{align*}
I_1&\leq C_{\kappa, c} \int_\Omega\int_{B(x, \frac{1}{2}d(x, \partial\Omega))}\phi\left(\frac{|u(x)-u(y)|}{{\frac{\lambda}{4C_1}}{|x-y|}^s}\right)\frac{dydx}{{|x-y|}^n}\\&
\leq \int_\Omega\int_{B(x, \frac{1}{2}d(x, \partial\Omega))}\phi\left(\frac{{\tilde{C}}|u(x)-u(y)|}{{{\lambda}}{|x-y|}^s}\right)\frac{dydx}{{|x-y|}^n}.
\end{align*}
For $I_2$, for each $Q\in \mathcal{F}$,  by (ii) $\forall Q\ne Q_0$,  we have $Q=Q_N$,
and
\begin{align*}
|u_{Q}- u_{Q_0}|&\leq\sum_{j=0}^{N-1} |u_{Q_j}-u_{Q_{j+1}}|\\&\leq\sum_{j=0}^{N-1} \left(|u_{Q_j}-u_{{Q_{j+1}}\cap{Q_j}}|+ |u_{Q_{j+1}}-u_{{Q_{j+1}}\cap{Q_j}}|\right).
\end{align*}
For adjacent cubes $Q_j, Q_{j+1}$,  one has

$$|{Q_j}- {Q_{j+1}}|\geq C_{\kappa, c}^{-1}\min\{|Q_j|, |Q_{j+1}|\}, $$
$$ C_{\kappa, c}^{-1}|Q_{j+1}|\leq|Q_j|\leq C_{\kappa, c}|Q_{j+1}|.$$
So
\begin{align*}
|u_{Q_j}-u_{{Q_{j+1}}\cap{Q_j}}|&\leq\frac{1}{|{Q_{j+1}}\cap{Q_j}|}\int_{Q_{j+1}\cap{Q_j}}|u(v)-u_{Q_j}|dv\\&
\leq\frac{{C_{\kappa, c}}^2}{|Q_j|}\int_{Q_j}|u(v)-u_{Q_j}|dv.
\end{align*}
Similarly,
$$|u_{Q_{j+1}}-u_{{Q_{j+1}}\cap{Q_j}}|\leq\frac{{C_{\kappa, c}}^2}{|Q_{j+1}|}\int_{Q_{j+1}}|u(v)-u_{Q_{j+1}}|dv.$$
As a result,  we get
$$|u_{Q}- u_{Q_0}|\leq2{C_{\kappa, c}}^2\sum_{j=0}^{N}\fint_{Q_j}|u(v)-u_{Q_{j}}|dv.$$
For each $Q_j$,   by the convexity of $\phi_\frac{n}{s}$ and Jenson inequality,
\begin{align*}
\fint_{Q_j}\frac{|u(v)-u_{Q_{j}}|}{\lambda}dv&={\phi_\frac{n}{s}}^{-1}\circ\phi_\frac{n}{s}\left(\fint_{Q_j}
\frac{|u(v)-u_{Q_{j}}|}{\lambda}dv\right)\\&
\leq{\phi_\frac{n}{s}}^{-1}\left(\fint_{Q_j}\phi_\frac{n}{s}
\left(\frac{|u(v)-u_{Q_{j}}|}{\lambda}\right)dv\right).
\end{align*}
By the inequality \eqref{eq2.1}, \eqref{gjs} and $\forall v, w \in Q_j,  |v-w|\leq\frac{1}{2}d(v, \partial\Omega), $
\begin{align*}
\int_{Q_j}\phi_\frac{n}{s}\left(\frac{|u(v)-u_{Q_{j}}|}{\lambda}\right)dv&\leq\int_{Q_j}\int_{Q_j}\phi
\left(\frac{|u(v)-u(w)|}{{\frac{\lambda}{C_1}}{|v-w|}^s}\right)\frac{dvdw}{{|v-w|}^n}\\&
\leq \int_{Q_j}\int_{B(v, \frac{1}{2}d(v, \partial\Omega))}\phi\left(\frac{|u(v)-u(w)|}{{\frac{\lambda}{C_1}}{|v-w|}^s}\right)\frac{dwdv}{{|v-w|}^n}
\\&:=\int_{Q_j}f(v)dv.
\end{align*}
Hence $$\fint_{Q_j}\frac{|u(v)-u_{Q_{j}}|}{\lambda}dv\leq{\phi_\frac{n}{s}}^{-1}\left(\fint_{Q_j}f(v)dv\right),$$
and
$$\frac{4|u_{Q}- u_{Q_0}|}{\lambda}\leq8{C_{\kappa, c}}^2\sum_{j=0}^{N}{\phi_\frac{n}{s}}^{-1}\left(\fint_{Q_j}f(v)dv\right).$$
By Lemma \ref{q1},
$$\phi_\frac{n}{s}\left(8{C_{\kappa, c}}^2\sum_{j=0}^{N}{\phi_\frac{n}{s}}^{-1}\left(\fint_{Q_j}f(v)dv\right)\right)\leq C\phi_\frac{n}{s}\left(\sum_{j=0}^{N}{\phi_\frac{n}{s}}^{-1}\left(\fint_{Q_j}f(v)dv\right)\right).$$
Applying $Q=Q_N\subset C_{\kappa, c}Q_j$ given in (ii),
$$|Q|\phi_\frac{n}{s}\left(\sum_{j=0}^{N}{\phi_\frac{n}{s}}^{-1}\left(\fint_{Q_j}f(v)dv\right)\right)
\leq\int_Q\phi_\frac{n}{s}\left(\sum_{P\in \mathcal{F}}{\phi_\frac{n}{s}}^{-1}\left(\fint_{P}f(v)dv\right)\chi_{C_{\kappa, c}P}\right)(x)dx.$$
Using the $\sum _{Q\in\mathcal{F}}\chi_{ Q}\leq\sum _{Q\in\mathcal{F}}\chi_{\kappa Q}  \leq  C_{\kappa, c} \chi_{\Omega}$ in (i) above,
\begin{align*}
I_2&\leq C\sum_{Q\in \mathcal{F}}\int_Q\phi_\frac{n}{s}\left(\sum_{P\in \mathcal{F}}{\phi_\frac{n}{s}}^{-1}\left(\fint_{P}f(v)dv\right)\chi_{C_{\kappa, c}P}\right)(x)dx\\
&\leq\tilde{C}\int_{\Omega}\phi_\frac{n}{s}\left(\sum_{P\in \mathcal{F}}{\phi_\frac{n}{s}}^{-1}\left(\fint_{P}f(v)dv\right)\chi_{C_{\kappa, c}P}\right)(x)dx.
\end{align*}
By Lemma \ref{le2.5},
$$I_2\leq C\int_{\Omega}\phi_\frac{n}{s}\left(\sum_{P\in \mathcal{F}}\left\{\mathcal{M}\left[\left({\phi_\frac{n}{s}}^{-1}(\fint_{P}f(v)dv)\right)^\frac{1}{2}
\chi_{P}\right]\right\}^2\right)(x)dx.$$
By Lemma \ref{le2.3},  $\phi_\frac{n}{s} \in \Delta_2\cap\nabla_2$.Let $\psi(t):=\phi_\frac{n}{s}(t^2)$,  then $\psi\in \Delta_2\cap\nabla_2$. Applying Lemma \ref{le2.4} to $q=2$ and $\psi$,  we obtain
$$I_2\leq C\int_{\Omega}\phi_\frac{n}{s}\left(\sum_{P\in \mathcal{F}}\left({\phi_\frac{n}{s}}^{-1}(\fint_{P}f(v)dv)\right)\chi_{P}\right)(x)dx.$$
Let $a_P=\fint_P f(v)dv$.   For each $x\in\Omega$,  we have
 \begin{align*}
\phi_\frac{n}{s}\left(\sum_{P\in \mathcal{F}}\left({\phi_\frac{n}{s}}^{-1}(a_P)\right)\chi_{P}(x)\right)
&=\phi_\frac{n}{s}\left(\frac{\sum_{P\in\mathcal{F}}\chi_{P}(x)}{\sum_{P\in\mathcal{F}}\chi_{P}(x)}\sum_{P\in \mathcal{F}}\left({\phi_\frac{n}{s}}^{-1}(a_P)\right)\chi_{P}(x)\right)\\&
\leq\phi_\frac{n}{s}\left(\frac{C_{\kappa, c}}{\sum_{P\in\mathcal{F}}\chi_{P}(x)}\sum_{P\in \mathcal{F}}\left({\phi_\frac{n}{s}}^{-1}(a_P)\right)\chi_{P}(x)\right)\\&
\leq\sum_{P\in\mathcal{F}}\frac{\chi_P(x)}{\sum_{P\in\mathcal{F}}\chi_{P}(x)}\phi_\frac{n}{s}(C_{\kappa, c}{\phi_\frac{n}{s}}^{-1}(a_P))\\&
\leq\sum_{P\in\mathcal{F}}\chi_P(x)\tilde{C}a_P.
\end{align*}
So
\begin{align*}
I_2&\leq C\int_\Omega\sum_{P\in\mathcal{F}}a_P\chi_P(x)dx\\&
\leq C\sum_{P\in\mathcal{F}}a_P|P|
=    C\sum_{P\in\mathcal{F}}\int_P f(v)dv\\&
\leq C(n, C_{\kappa, c}, K_\phi)\int_\Omega\int_{B(v, \frac{1}{2}d(v, \partial\Omega))}\phi\left(\frac{C|u(v)-u(y)|}{\lambda{|u-w|}^s}\right)\frac{dwdv}{|u-w|^n}.
\end{align*}
In the end,  we obtain
$$I\leq\int_\Omega\int_{B(v, \frac{1}{2}d(v, \partial\Omega))}\phi\left(\frac{C|u(v)-u(y)|}{\lambda{|u-w|}^s}\right)\frac{dwdv}{|u-w|^n}, $$
where $C=C(n,  s, K_\phi)$,  $\lambda>4C_1\|u\|_{{\dot {V}_\ast}^{s, \phi}(\Omega)}$and $C\geq 4C_1$.
Let $\lambda>C\|u\|_{{\dot {V}_\ast}^{s, \phi}(\Omega)}, ~$we have $I\leq1$.
\end{proof}

\section{Proof of Themrem \ref{th1}(ii) }
To prove Theorem \ref{th1} (ii),  the most important method is getting the fact which Lemma \ref{le421} expressed. We first need to choose a special test function to estimate the relationship between its norms and its radius.

Let $z\in \Omega,  ~d(z, \partial\Omega)\leq m <\diam\Omega$. Denote $\Omega_{z, m}$ by a component of $\Omega\setminus{\overline {B_{\Omega}(z, m)}}$.
For $t>r\geq m$ with $\Omega_{z, m}\neq\varnothing,  $
define $u_{z, r, t}$ in $\Omega$ as
\begin{align}\label{3.1}
u_{z, r, t}(y)=\left\{\begin{array}{ll}
0& y\in \Omega\setminus [\Omega_{z, m} \setminus B_\Omega(z, r)]\\
\frac{|y-z|-r}{t-r} & y\in \Omega_{z, m}\cap [B(z, t)\setminus B(z, r)], \\
1& y\in\Omega_{z, m} \setminus B_\Omega(z, t),
\end{array}\right.
\end{align}
where $B_\Omega(z, t)=B(z, t)\cap\Omega$.

It's not difficult to know the following property.
\begin{lem}
$u_{z, r, t}$ is Lipschitz with the Lipschitz constant  $\frac{1}{t-r}$.
\end{lem}
\begin{proof} We spilt into three cases to prove it.

Case 1. For $x \in\Omega\setminus [\Omega_{z, m} \setminus B_\Omega(z, r)]$,
it means that $u_{z, r, t}(x)=0$.
Since $u_{z, r, t}(y)=u_{z, r, t}(x)=0$ when $y \in\Omega\setminus [\Omega_{z, m} \setminus B_\Omega(z, r)]$,
we only need to consider $y \in \Omega_{z, m}\cap [B(z, t)\setminus B(z, r)]$ or $y \in \Omega_{z, m} \setminus B_\Omega(z, t)$.
If $y \in \Omega_{z, m}\cap [B(z, t)\setminus B(z, r)]$,  we know $|x-z|\leq r$. Hence
\begin{align*}
|u_{z, r, t}(x)-u_{z, r, t}(y)|=\frac{|y-z|-r}{t-r}\leq\frac{|y-z|-|x-z|}{t-r}\leq\frac{|x-y|}{t-r}.
\end{align*}
If  $y \in \Omega_{z, m} \setminus B_\Omega(z, t)$,  we get $|x-y|\geq t-r$. Therefore,
\begin{align*}
|u_{z, r, t}(x)-u_{z, r, t}(y)|=1 \leq\frac{|x-y|}{t-r}.
\end{align*}

Case 2. For $x \in \Omega_{z, m}\cap [B(z, t)\setminus B(z, r)]$,
then $u_{z, r, t}(x)=\frac{|x-z|-r}{t-r}$.
If $y \in \Omega_{z, m}\cap [B(z, t)\setminus B(z, r)]$ with $u_{z, r, t}(y)=\frac{|y-z|-r}{t-r}$,
\begin{align*}
|u_{z, r, t}(x)-u_{z, r, t}(y)|=\left|\frac{|x-z|-r}{t-r}-\frac{|y-z|-r}{t-r}\right|\leq\frac{|x-z|-|y-z|}{t-r}\leq\frac{|x-y|}{t-r}.
\end{align*}
If $y \in\Omega\setminus [\Omega_{z, m} \setminus B_\Omega(z, r)]$ with $u_{z, r, t}(y)=0$,
we have $|y-z|\leq r$.
Then
\begin{align*}
|u_{z, r, t}(x)-u_{z, r, t}(y)|=\frac{|x-z|-r}{t-r}\leq\frac{|y-z|-|x-z|}{t-r}\leq\frac{|x-y|}{t-r}.
\end{align*}
If $y \in \Omega_{z, m} \setminus B_\Omega(z, t)$ with $u_{z, r, t}(y)=1$,
then $|y-z|\geq t$.
Together with $|x-z|\leq t$,  we have
\begin{align*}
|u_{z, r, t}(x)-u_{z, r, t}(y)|&=\left|\frac{|x-z|-r}{t-r}-1\right|=\left|\frac{|x-z|-t}{t-r}\right|
=\frac{t-|x-z|}{t-r}\\
&\leq\frac{|y-z|-|x-z|}{t-r}\leq\frac{|x-y|}{t-r}.
\end{align*}

Case 3. For $x \in \Omega_{z, m} \setminus B_\Omega(z, t)$,
then $u_{z, r, t}(x)=1$.
Since $u_{z, r, t}(y)=u_{z, r, t}(x)=1$ when $y \in\Omega_{z, m} \setminus B_\Omega(z, t)$,
we only need to consider $y \in \Omega\setminus [\Omega_{z, m} \setminus B_\Omega(z, r)]$ or $y \in \Omega_{z, m}\cap [B(z, t)\setminus B(z, r)]$.
If $y \in \Omega\setminus [\Omega_{z, m} \setminus B_\Omega(z, r)]$ with $u_{z, r, t}(y)=0$,
together with $|x-y|\geq t-r$,  we know
\begin{align*}
|u_{z, r, t}(x)-u_{z, r, t}(y)|=1 \leq\frac{|x-y|}{t-r}.
\end{align*}
If  $y \in \Omega_{z, m}\cap [B(z, t)\setminus B(z, r)]$ with $u_{z, r, t}(y)=\frac{|y-z|-r}{t-r}$,
then $|y-z|\leq t$.
Moreover,  $|x-z|\geq t$. Hence
\begin{align*}
|u_{z, r, t}(x)-u_{z, r, t}(y)|&=\left|1-\frac{|y-z|-r}{t-r}\right|\\
&\leq\frac{|x-z|-|y-z|}{t-r}\leq\frac{|x-y|}{t-r}.
\end{align*}

Combining above cases,  we know $u_{z, r, t}$ is Lipschitz with the Lipschitz constant  $\frac{1}{t-r}$.
\end{proof}

Next we provide an estimation of the test function.
\begin{lem}\label{le3.x}
Let $s\in(0, 1)$ and $\phi$ be a Young function. For any bounded domain $\Omega\subset\mathbb{R}^n$ and $z\in \Omega$ with $d(z, \partial\Omega)\leq m<\diam\Omega$. For $t>r\geq m$,  we have $u_{z, r, t}\in{\dot {V}_\ast}^{s, \phi}(\Omega)$  with
$$\|u_{z, r, t}\|_{{\dot {V}_\ast}^{s, \phi}(\Omega)}\leq C\left({\phi}^{-1}\left(\frac{1}{|\Omega_{z, m}\setminus B(z, r)|}\right)\right)^{-1}\frac{1}{(t-r)^s}, $$
where $C=C(n, s, C_\phi)\geq 1$.
\end{lem}
\begin{proof}
For any $x\in\Omega$ and $y\in B(x, \frac{1}{2}d(x, \partial\Omega))\subset\Omega, |u_{z, r, t}(x)-u_{z, r, t}(y)|\neq 0$ means that either $x$ or $y$ in $\Omega_{z, m}\setminus B(z, r)$.
\begin{align*}
H:&=\int_\Omega\int_{|x-y|<\frac{1}{2}d(x, \partial\Omega)}
\phi\left(\frac{|u_{z, r, t}(x)-u_{z, r, t}(y)|}{\lambda|x-y|^s}\right)\frac{dydx}{|x-y|^n}\\&
\leq 2\int_{\Omega_{z, m}\setminus B(z, r)}\int_\Omega\phi\left(\frac{|u_{z, r, t}(x)-u_{z, r, t}(y)|}
{\lambda|x-y|^s}\right)\frac{dydx}{|x-y|^n}\\&
\leq2\int_{\Omega_{z, m}\setminus B(z, r)}\int_{B(x, t-r)}\phi\left(\frac{|x-y|^{1-s}}{\lambda(t-r)}\right)\frac{dydx}{|x-y|^n}\\&
+2\int_{\Omega_{z, m}\setminus B(z, r)}\int_{\mathbb{R}^n\setminus B(x, t-r)}\phi\left(\frac{1}{\lambda|x-y|^s}\right)\frac{dydx}{|x-y|^n}\\&
:=2H_1+2H_2.
\end{align*}
Using change of variable and \eqref{eq1},  we have
\begin{align*}
H_1&=\int_{\Omega_{z, m}\setminus B(z, r)}\int_0^{t-r}n\omega_n\phi\left(\frac{\rho^{1-s}}{\lambda(t-r)}\right)\frac{d\rho}{\rho}dx\\&
=\int_{\Omega_{z, m}\setminus B(z, r)}\int_0^{\frac{1}{\lambda(t-r)^{s}}}n\omega_n\frac{1}{1-s}\phi(\mu)\frac{d\mu}{\mu}dx\\&
\leq\int_{\Omega_{z, m}\setminus B(z, r)}\frac{C_\phi n\omega_n}{1-s}\phi\left(\frac{1}{\lambda(t-r)^s}\right)dx\\&
=\frac{C_\phi n\omega_n}{1-s}\phi\left(\frac{1}{\lambda(t-r)^s}\right)|\Omega_{z, m}\setminus B(z, r)|,
\end{align*}
and
\begin{align*}
H_2&=\int_{\Omega_{z, m}\setminus B(z, r)}\int_{t-r}^\infty n\omega_n\phi\left(\frac{1}{\lambda\rho^{s}}\right)\frac{d\rho}{\rho}dx\\&
=\int_{\Omega_{z, m}\setminus B(z, r)}\int_0^{\frac{1}{\lambda(t-r)^{s}}}n\omega_n\frac{1}{s}\phi(\mu)\frac{d\mu}{\mu}dx\\&
\leq\int_{\Omega_{z, m}\setminus B(z, r)}\frac{C_\phi n\omega_n}{s}\phi\left(\frac{1}{\lambda(t-r)^s}\right)dx\\&
=\frac{C_\phi n\omega_n}{s}\phi\left(\frac{1}{\lambda(t-r)^s}\right)|\Omega_{z, m}\setminus B(z, r)|.
\end{align*}
Let $\lambda=M\left({\phi}^{-1}\left(\frac{1}{|\Omega_{z, m}\setminus B(z, r)|}\right)\right)^{-1}\frac{1}{(t-r)^s}$,  where $M\geq \max\left\{\frac{4C_\phi n\omega_n}{1-s}, \frac{4C_\phi n\omega_n}{s}, 1\right\}$. We have $H_1\leq\frac{1}{4}, H_2\leq\frac{1}{4}$, hence $H\leq1$.
As a result $$\|u_{z, r, t}\|_{{\dot {V}_\ast}^{s, \phi}(\Omega)}\leq C\left({\phi}^{-1}\left(\frac{1}{|\Omega_{z, m}\setminus B(z, r)|}\right)\right)^{-1}\frac{1}{(t-r)^s}.$$
\end{proof}

For $x_0, z\in \Omega, $ let $r>0$ such that $d(z, \partial\Omega)<r<|x_0-z|$.   Define \
$$\omega_{x_0, z, r}(y):=\frac{1}{r}\inf_{\gamma(x_0, y)}{l(\gamma\cap B(z, r))},  \forall y \in \Omega, $$
where the infimum is taken over all reactiable curves $\gamma$ joining $x_0$ and $y$.

\begin{lem}\label{le3.3}
 $s\in(0, 1)$ and $\phi$ be a Young function. For any bounded domain $\Omega\subset\mathbb{R}^n$ and $x_0, z\in \Omega$ and $r>0$ with $d(z, \partial\Omega)\leq r< |x_0-z|, $ we have $\omega_{x_0, z, r}\in{\dot {V}_\ast}^{s, \phi}(\Omega)$  with
$$\|\omega_{x_0, z, r}\|_{{\dot {V}_\ast}^{s, \phi}(\Omega)}\leq C\left({\phi}^{-1}\left(\frac{1}{r^n}\right)\right)^{-1}\frac{1}{r^s}, $$
where $C=C(n, s, C_\phi)\geq 1$.
\end{lem}
\begin{proof}
For $x\in\Omega\setminus B(z, 6r),  y\in B(x, \frac{1}{2}d(x, \partial\Omega)), $ we have
$$d(x, \partial\Omega)\leq|x-z|+d(z, \partial\Omega)\leq|x-z|+r, $$and
\begin{align*}
|y-z|&\geq|x-z|-|y-x|\\&\geq|x-z|-\frac{1}{2}(|x-z|+r)\\&=\frac{1}{2}|x-z|-\frac{r}{2}
\\&\geq3r-\frac{r}{2}\geq2r.
\end{align*}
So $B(x, \partial\Omega))\cap B(z, 2r)=\varnothing$.  ~Let $\gamma_{x, y}$ be the segment joining $x, y$ containing in $B(x, \frac{1}{2}d(x, \partial\Omega)), $ then  $\gamma_{x, y}\subset\Omega\setminus\ B(z, r)$.
For any $\gamma(x_0, x), \gamma(x_0, x)\cup\gamma_{x, y}$ is a curve joining $x_0$ and $y$,  with
$$l((\gamma(x_0, x)\cup\gamma_{x, y})\cap B(z, r))=l(\gamma(x_0, x)\cap B(z, r)).$$
Hence $\omega_{x_0, z, r}(y)\leq\omega_{x_0, z, r}(x)$.

Similarity  $\omega_{x_0, z, r}(x)\leq\omega_{x_0, z, r}(y)$.
So for any $x\in \Omega\setminus B(z, 6r) , y\in B(x, \frac{1}{2}d(x, \partial\Omega))$,  we have $\omega_{x_0, z, r}(x)=\omega_{x_0, z, r}(y)$.

For any $x\in \Omega, ~|x-y|<\frac{1}{2}d(x, \partial\Omega), $ we know $l(\gamma_{x, y}\cap B(z, r))\leq|x-y|$.   Since $\gamma(x_0, x)\cup\gamma_{x, y}$ is a curve joining $x_0$ and $y$,  we get $$\omega_{x_0, z, r}(y)\leq\omega_{x_0, z, r}(x)+\frac{1}{r}|x-y|.$$
Similarity $\omega_{x_0, z, r}(x)\leq\omega_{x_0, z, r}(y)+\frac{1}{r}|x-y|$. So $|\omega_{x_0, z, r}(y)-\omega_{x_0, z, r}(x)|\leq\frac{1}{r}|x-y|$.

For $x\in\Omega\cap B(z, 6r), $ we have $d(x, \partial\Omega)\leq6r+d(z, \partial\Omega)<8r$.
\begin{align*}
H:&=\int_\Omega\int_{|x-y|<\frac{1}{2}d(x, \partial\Omega)}\phi\left(\frac{|\omega_{x_0, z, r}(x)-\omega_{x_0, z, r}(y)|}{\lambda|x-y|^s}\right)\frac{dydx}{|x-y|^n}\\&
=\int_{\Omega\cap B(z, 6r)}\int_{|x-y|<\frac{1}{2}d(x, \partial\Omega)}\phi\left(\frac{|\omega_{x_0, z, r}(x)-\omega_{x_0, z, r}(y)|}{\lambda|x-y|^s}\right)\frac{dydx}{|x-y|^n}\\&
\leq\int_{\Omega\cap B(z, 6r)}\int_0^{4r}n\omega_n\phi\left(\frac{\rho^{1-s}}{\lambda}\right)\frac{d\rho}{\rho}dx\\&
\leq\int_{\Omega\cap B(z, 6r)}\frac{C_\phi n\omega_n}{1-s}\phi\left(\frac{4^{1-s}}{\lambda r^s}\right)dx\\&
\leq\frac{C_\phi n\omega_n^2}{1-s}\phi\left(\frac{4^{1-s}}{\lambda r^s}\right)(6r)^n
\end{align*}
Let $\lambda=M\left(\phi^{-1}\left(\frac{1}{r^n})\right)\right)^{-1}\frac{1}{r^s} $,  where $M>\max\left\{\frac{C_\phi n\omega_n^24^{1-s}}{1-s}6^n, 4^{1-s}\right\} $, then $H\leq1$.
So $$\|\omega_{x_0, z, r}\|_{{\dot {V}_\ast}^{s, \phi}(\Omega)}\leq C\left({\phi}^{-1}\left(\frac{1}{r^n}\right)\right)^{-1}\frac{1}{r^s}.$$
\end{proof}
\begin{lem}\label{ls48}
Let $s\in(0, 1)$ and $\phi \in \Delta_2$ be a Young function satisfying $K_{\phi}<2^{\frac{n}{s}}$ in $\eqref{da2}$,  a bounded domain $\Omega\subset\mathbb{R}^n$ supports the $(\phi_\frac{n}{s}, \phi)$- Poincar\'e inequality $\eqref{ueq1}$. Fix a point $x_0$ so that $r_0:=\max\{d(x, \partial\Omega):x\in \Omega\}=d(x_0, \partial\Omega)$.  Assume that $x,  x_0\in \Omega\setminus{\overline {B(z, r)}}$ for some $z\in \Omega$ and $r\in(0, 2\diam\Omega)$,  there exists a positive constant $b_0$ that $x, x_0$ are contained in the same component of $\Omega\setminus{\overline {B(z, b_0r)}}$.
\end{lem}
\begin{proof}
Let $b_{x, z, r} :=\sup\left\{c\in(0, 1]:x, x_0 ~in ~the ~same~ component ~of ~\Omega\setminus{\overline {B(z, cr)}}\right\}$.
We need prove that $b_{x, z, r}$ has the positive low bound independent of $x, z, r$, that is
$$b=\inf\left\{b_{x, z, r}:\exists z\in \Omega ,  r\in(0, 2\diam\Omega) ~such ~that~ x, x_0 \in \Omega\setminus{\overline {B(z, r)}}\right\}>0.$$
then let $b_0=\frac{b}{2}$,  we get the conclusion.
Because it is a infimum problem,  we may assume $b_{x, z, r}\leq\frac{1}{10}$.

We want to prove $$\frac{r}{C}(\frac{1}{2}-2b_{x, z, r})\leq|\Omega_x|^\frac{1}{n}\leq 2Cb_{x, z, r}r,  ~C\geq 1.$$
then $$b_{x, z, r}\geq\frac{1}{4(C^2+1)},$$ so $b>0$.

First for fixed $x, z, r$,  we have $b_{x, z, r}>0$.
By $z\in \Omega $,  then there existing $ 0<\delta<1$ such that $B(z, \delta r)\subset\Omega, $ and $x_0 \notin B(z, \delta r)$.
For $h=\frac{\delta}{2}$,  and a curve $\gamma(x, x_0)$ if
$$\gamma(x, x_0)\cap\overline{B(z, hr)}=\varnothing, $$ then $x, x_0$ are contained in the same component of $\Omega\setminus{\overline {B(z, hr)}}$.

If $\gamma(x, x_0)\cap\overline{B(z, hr)}\neq\varnothing$,
denote $t_0:=\inf\left\{ t\in [0, 1]:\gamma(x, x_0)(t)\in\partial B(z, \delta r)\right\}$,\\
$t_1:=\sup\left\{ t\in [0, 1]:\gamma(x, x_0)(t)\in\partial B(z, \delta r)\right\}$
and $A:=\gamma(x, x_0)(t_0), B:=\gamma(x, x_0)(t_1)$.
Then we have
\begin{align*}
\tilde{\gamma}=\gamma(x, x_0)|_{t\in(0, t_0)}\cup\overset{\frown}{AB}\cup\gamma(x, x_0)|_{t\in(t_0, 1)}\subset
\Omega\setminus{\overline {B(z, hr)}}.
\end{align*}
and  $x, x_0$ are contained in the same component of $\Omega\setminus{\overline {B(z, hr)}}$.
So $b_{x, z, r}\geq h>0$.

Set $c_0=2b_{x, z, r}\leq\frac{1}{5}, $then $x_0\notin{\overline {B(z, c_0r)}}$.  Denote by $\Omega_{x_0}$ the component of $\Omega\setminus{\overline {B(z, c_0r)}}$ containing $x_0$.
By $b_{x, z, r}<\frac{2}{3}c_0<1, $ we have $x, x_0$ are not contained in the same component of $\Omega\setminus{\overline {B(z, \frac{2}{3}c_0r)}}$. Now we prove that $B(z, c_0r)\cap\partial\Omega\neq\varnothing$.
If not,  by $z\in \Omega, $ we have
$$B(z, \frac{2}{3}c_0)\subset B(z, c_0r)\subset\Omega.$$
From the above discussion,  we get $x, x_0$ are  contained in the same component of $\Omega\setminus{\overline {B(z, \frac{2}{3}c_0r)}}, $ and we get contradiction. So $B(z, c_0r)\cap\partial\Omega\neq\varnothing$.
Then $$r_0=d(x_0, \partial\Omega)\leq\underset{y\in B(z, c_0r)}{\max}|x_0-y|\leq r+c_0r+d(x_0, B(z, r))\leq\frac{6}{5}r+d(x_0, B(z, r)).$$
$$d(x_0, B(z, c_0r))\geq|x_0-z|-\frac{r}{5}=d(x_0, B(z, r))+\frac{4}{5}r.$$
So $d(x_0, B(z, c_0r))\geq\frac{r_0}{2}, ~$and
\begin{equation}\label{eq3.1}
B(x_0, \frac{r_0}{2})\subset\Omega_{x_0}\subset\Omega\setminus\Omega_{x}.
\end{equation}

Define
$$\omega(y):=\frac{1}{c_0r}\underset{\gamma(x_0, y)}{inf}l(\gamma\cap B(z, c_0r)), ~\forall y\in\Omega.$$
Since $B(z, c_0r)\cap\partial\Omega\neq\varnothing$ and $x_0\notin{\overline {B(z, c_0r)}},$ we have
$d(z, \partial\Omega)<c_0r<|x_0-z|$.
By Lemma \ref{le3.3},  we know $$\|\omega\|_{{\dot {V}_\ast}^{s, \phi}(\Omega)}\leq C\left({\phi}^{-1}\left(\frac{1}{(c_0r)^n}\right)\right)^{-1}\frac{1}{(c_0r)^s}, $$
By the $(\phi_\frac{n}{s}, \phi)$- Poincare inequality \eqref{ueq1} ,
$$\|\omega-\omega_\Omega\|_{L^{\phi_\frac{n}{s}}(\Omega)}\leq C\|\omega\|_{{\dot {V}}^{s, \phi}(\Omega)}\leq C\left({\phi}^{-1}\left(\frac{1}{(c_0r)^n}\right)\right)^{-1}\frac{1}{(c_0r)^s}.$$
On the other hand, by \eqref{eq3.1},  $y\in B(x_0, \frac{1}{2}r_0),  \omega(y)=0$.
Since $\Omega$ is bounded,  $r_0>0,$ we have $\frac{|\diam\Omega|}{r_0^n}\leq C$.
Using the convexity of $\phi_{\frac{n}{s}}, $
\begin{align*}
\int_\Omega\phi_{\frac{n}{s}}\left(\frac{|\omega(x)|}{\lambda}\right)dx\leq\frac{1}{2}\int_\Omega
\phi_{\frac{n}{s}}\left(\frac{|\omega(x)-\omega_\Omega|}{\lambda}\right)dx
+\frac{|\Omega|}{2}\phi_{\frac{n}{s}}\left(\frac{|\omega_{B(x_0, \frac{1}{2}r_0)}-\omega_\Omega|}{\lambda}\right).
\end{align*}
By the Jensen inequality,
\begin{align*}
|\Omega|\phi_{\frac{n}{s}}\left(\frac{|\omega_{B(x_0, \frac{1}{2}r_0)}-\omega_\Omega|}{\lambda}\right)&
\leq|\Omega|\fint_{B(x_0, \frac{1}{2}r_0)}\phi_{\frac{n}{s}}\left(\frac{|\omega(x)-\omega_\Omega|}{\lambda}\right)dx\\&
\leq\frac{|\Omega|}{|B(x_0, \frac{1}{2}r_0)|}\int_\Omega\phi_{\frac{n}{s}}\left(\frac{|\omega(x)-\omega_\Omega|}{\lambda}\right)dx\\&
\leq2^nC^n\int_\Omega\phi_{\frac{n}{s}}\left(\frac{|\omega(x)-\omega_\Omega|}{\lambda}\right)dx.
\end{align*}
As a result,
$$\int_\Omega\phi_{\frac{n}{s}}\left(\frac{|\omega(x)|}{\lambda}\right)dx\leq C\int_\Omega\phi_{\frac{n}{s}}\left(\frac{|\omega(x)-\omega_\Omega|}{\lambda}\right)dx, $$ and
\begin{equation}\label{t23}
\|\omega\|_{L^{\phi_\frac{n}{s}}(\Omega)}\leq C\|\omega-\omega_\Omega\|_{L^{\phi_\frac{n}{s}}(\Omega)}.
\end{equation}

Since $\forall y\in \Omega_x, ~\omega(y)\geq1, $
we have $$\int_\Omega\phi_{\frac{n}{s}}\left(\frac{|\omega(x)|}{\lambda}\right)dx\geq\phi_{\frac{n}{s}}
\left(\frac{1}{\lambda}\right)|\Omega_x|, $$
and $$\|\omega\|_{L^{\phi_\frac{n}{s}}(\Omega)}\geq \left({\phi_{\frac{n}{s}}}^{-1}\left(\frac{1}{|\Omega_x|}\right)\right)^{-1}.$$
So $$C\phi^{-1}\left(\frac{1}{(c_0r)^n}\right)(c_0r)^s
\leq{\phi_{\frac{n}{s}}}^{-1}\left(\frac{1}{|\Omega_x|}\right).$$
By \eqref{te1}, $$\frac{H(A)}{A}\leq C\frac{1}{{\phi(A)}^\frac{s}{n}}.$$
Let $$A=\phi^{-1}\left(\frac{1}{(c_0r)^n}\right), $$
we have
$$\frac{{\phi_\frac{n}{s}}^{-1}\left(\frac{1}{(c_0r)^n}\right)}{\phi^{-1}\left(\frac{1}{(c_0r)^n}\right)}\leq C(c_0r)^s.$$
So $${\phi_\frac{n}{s}}^{-1}\left(\frac{1}{(c_0r)^n}\right)\leq C{\phi_{\frac{n}{s}}}^{-1}\left(\frac{1}{|\Omega_x|}\right).$$
By Lemma \ref{le2.3}, $\phi_\frac{n}{s}\in \Delta_2, $ and Lemma \ref{q1},
we have $$\frac{1}{(c_0r)^n}\leq C\frac{1}{|\Omega_x|}, $$
and \begin{equation}\label{233}
{|\Omega_x|}^\frac{1}{n}\leq C(c_0r).
\end{equation}
Let $c_j>c_{j-i}$ for $j\geq1$ such that $$|\Omega_x\setminus{B(z, c_jr)}|=\frac{1}{2}|\Omega_x\setminus B(z, c_{j-1}r)|=2^{-j}|\Omega_x|.$$
For $j\geq0$ with $\Omega_x\setminus \overline{B(z, c_jr)}\neq\varnothing, $
define $v_j$ in $\Omega$ as
\begin{align*}
v_j(y)=\left\{\begin{array}{ll}
0& y\in \Omega\setminus [\Omega_x \setminus B_\Omega(z, c_{j+1}r)]\\
\frac{|y-z|-c_jr}{c_{j+1}r-c_jr} & y\in \Omega_x\cap [B(z, c_jr)\setminus B(z, c_{j+1}r)], \\
1& y\in\Omega_x \setminus B_\Omega(z, c_jr),
\end{array}\right.
\end{align*}
Let $\Omega_{z, x}=\Omega_x, r=c_jr$ and $t=c_{j+1}r$, then $v_j(y)=u_{z, c_jr, c_{j+1}r}(y)$ where $u_{z, c_jr, c_{j+1}r}(y)$ is defined in \eqref{3.1}.
Applying Lemma \ref{le3.x},  we have
$$\|v_j\|_{{\dot {V}_\ast}^{s, \phi}(\Omega)}\leq C\left({\phi}^{-1}\left(\frac{1}{|\Omega_x\setminus B(z, c_jr)|}\right)\right)^{-1}\frac{1}{(c_{j+1}r-c_jr)^s}.$$
Applying \eqref{eq3.1}, we have $v_j(y)=0$ for $y\in B(x_0, \frac{1}{2}r_0)$.   Similarly to \eqref{t23},  we have
\begin{equation}
\|v_j\|_{L^{\phi_\frac{n}{s}}(\Omega)}\leq C\|v_j-{v_j}_\Omega\|_{L^{\phi_\frac{n}{s}}(\Omega)}.
\end{equation}
And $v_j(y)=1$ for $y\in \Omega_x\setminus {B_\Omega(z, c_jr)}$,  then we have
$$\|v_j\|_{L^{\phi_\frac{n}{s}}(\Omega)}\geq \left({\phi_{\frac{n}{s}}}^{-1}\left(\frac{1}{|\Omega_x\setminus {B_\Omega(z, c_jr)}|}\right)\right)^{-1}.$$
By the $(\phi_\frac{n}{s}, \phi)$-Poincar\'e inequality \eqref{ueq1} , we have
$${\phi_{\frac{n}{s}}}^{-1}\left(\frac{1}{|\Omega_x\setminus {B_\Omega(z, c_jr)}|}\right)\geq C{\phi}^{-1}\left(\frac{1}{|\Omega_x\setminus B(z, c_jr)|}\right)(c_{j+1}r-c_jr)^s.$$
By \eqref{te1}, $$\frac{H(A)}{A}\leq C\frac{1}{{\phi(A)}^\frac{s}{n}}.$$
and let $$A={\phi}^{-1}\left(\frac{1}{|\Omega_x\setminus B(z, c_jr)|}\right), $$then
$$(c_{j+1}r-c_jr)^s\leq C|\Omega_x\setminus B(z, c_jr)|^\frac{s}{n}.$$
So $c_{j+1}-c_jr\leq C|\Omega_x\setminus B(z, c_jr)|^\frac{1}{n}\leq C2^{-\frac{j}{n}}|\Omega_x|^\frac{1}{n}$.

Now we prove that $\sup\left\{c_j\right\}>1$.
If not, we have $\forall c_j\leq1$.
By $x\in\Omega\setminus\overline{B(x, r)},$ then $\exists \delta>0$ such that $$B(x, \delta)\subset\Omega\setminus\overline{B(x, r)}\subset\Omega\setminus\overline{B(x, c_0r)}.$$
By the connectivity of the $B(x, \delta)$,  we have $B(x, \delta)\subset\Omega_x$.
Then $$B(x, \delta)\subset\Omega_x\setminus\overline{B(x, r)}\subset\Omega_x\setminus B(x, c_jr), $$
and
$$0<|B(x, \delta)|\leq|\Omega_x\setminus\overline{B(x, r)}|\leq|\Omega_x\setminus B(x, c_jr)|=2^{-j}|\Omega_x|.$$
Letting $j\to \infty$, we get a contradiction,
and hence $\sup\left\{c_j\right\}>1$.
So there exists $c_j$ such that $c_j\geq\frac{1}{2}.$ Let $j_0=\inf\left\{j\geq1:c_j\leq\frac{1}{2}\right\}, $
then $$(\frac{1}{2}-c_0)r\leq(c_{j_0}-c_0)r=\sum_{j=0}^{j_0-1}(c_{j+1}-c_j)r\leq C\sum_{j=0}^{j_0-1}2^{-\frac{j}{n}}|\Omega_x|^\frac{1}{n}\leq 2C|\Omega_x|^\frac{1}{n}.$$
So $\frac{r}{C}(\frac{1}{2}-2b_{x, z, r})\leq|\Omega_x|^\frac{1}{n}$.
By the \eqref{233},  we have
$$\frac{r}{C}(\frac{1}{2}-2b_{x, z, r})\leq|\Omega_x|^\frac{1}{n}\leq C2b_{x, z, r}r,  \, C\geq 1.$$
Then $$b_{x, z, r}\geq\frac{1}{4(C^2+1)}, $$which implies $b>0$.
\end{proof}
\begin{lem}\label{le421}
Let $s\in(0, 1)$ and $\phi \in \Delta_2$ be a Young function satisfying $K_\phi<2^{\frac{n}{s}}$ in \eqref{da2}.  If a bounded domain $\Omega\subset\mathbb{R}^n$ supports the $(\phi_\frac{n}{s}, \phi)$- Poincar\'e inequality \eqref{ueq1}, then the $\Omega$ has the LLC(2) property,  that is,  there exists a constant $b\in(0, 1)$ such that for all $z\in \mathbb{R}^n$ and $r>0, $ any pair of point in $\Omega\setminus\overline{B(z, r)}$ can be joined in $\Omega\setminus\overline{B(z, br)}$.
\end{lem}

\begin{proof}
Fix $x_0$ so that $r_0:=\max(d(x, \partial\Omega):x\in \Omega)=d(x_0, \partial\Omega))$
and $b_0$ is the constant in Lemma \ref{ls48}.
Then we spilt into three cases to prove it.

Case 1. For $z\notin B\left(x_0, \frac{r_0}{8\diam\Omega}r\right)$, we consider the radius $r$.

If $r>\frac{16(\diam\Omega)^2}{r_0}$, then $\forall y \in\overline{ B\left(z, \frac{r_0}{16\diam\Omega}r\right)}$, we have
$$|y-x_0|\geq|z-x_0|-|z-y|\geq\frac{r_0}{16\diam\Omega}r>\diam\Omega.$$
By $\Omega\subset B(x_0, \diam\Omega)$, we get $\Omega\cap\overline{ B\left(z, \frac{r_0}{16\diam\Omega}r\right)}=\varnothing$. Here, any pair of point in $\Omega\setminus\overline{B(z, r)}$ can be joined in $\Omega\setminus\overline{B(z, \frac{r_0}{16\diam\Omega}r)}=\Omega$.

If $r\leq\frac{16(\diam\Omega)^2}{r_0}$ and $d(z, \partial\Omega)>\frac{b_0r_0}{32\diam\Omega}r$.
When $z\notin\Omega$, then any pair of point in $\Omega\setminus\overline{B(z, r)}$ can be joined in $\Omega\setminus\overline{B\left(z, \frac{b_0r_0}{32\diam\Omega}r\right)}=\Omega$.
When $z\in\Omega$, then $B\left(z, \frac{b_0r_0}{64\diam\Omega}r\right)\subset B\left(z, \frac{b_0r_0}{32\diam\Omega}r\right)\subset \Omega$.
Similar to the process of proving $b_{x, z, r}>0$ in Lemma \ref{ls48},  we know $\Omega\setminus\overline{B\left(z, \frac{b_0r_0}{64\diam\Omega}r\right)}$ is a connected set.
Here,  any pair of point in $\Omega\setminus\overline{B(z, r)}$ can be joined in $\Omega\setminus\overline{B\left(z, \frac{b_0r_0}{64\diam\Omega}r\right)}$.

If $r\leq\frac{16(\diam\Omega)^2}{r_0}$ and $d(z, \partial\Omega)\leq\frac{b_0r_0}{32\diam\Omega}r$.
Let $y\in B\left(z, \frac{b_0r_0}{16\diam\Omega}r\right)\cap\Omega$.
By $B\left(y, (1-\frac{b_0}{2})\frac{r_0}{8\diam\Omega}r\right)\subset B\left(z, \frac{r_0}{8\diam\Omega}r\right)\subset B(z, r)$, we know
$${\forall x\in \Omega\setminus \overline {B(z, r)}, x, x_0\in\Omega\setminus \overline {B\left(y, (1-\frac{b_0}{2})\frac{r_0}{8\diam\Omega}r\right)}.}$$
By Lemma \ref{ls48}, $x, x_0$ are in the same component of $\Omega\setminus{\overline {B\left(y, b_0(1-\frac{b_0}{2})\frac{r_0}{8\diam\Omega}r\right)}}$. By $$\forall w\in B\left(z, \frac{b_0(1-b_0)r_0}{16\diam\Omega}r\right), $$ we have
$$|w-y|\leq|w-z|+|z-y|<\frac{b_0(1-b_0)r_0}{16\diam\Omega}r+\frac{b_0r_0}{16\diam\Omega}r
=b_0\left(1-\frac{b_0}{2}\right)\frac{r_0}{8\diam\Omega}r.$$
Then $$B\left(z, \frac{b_0(1-b_0)r_0}{16\diam\Omega}r\right)\subset B\left(y, b_0\left(1-\frac{b_0}{2}\right)\frac{r_0}{8\diam\Omega}r\right), $$
and $\Omega\setminus{\overline { B\left(y, b_0\left(1-\frac{b_0}{2}\right)\frac{r_0}{8\diam\Omega}r\right)}}\subset\Omega\setminus\overline {B\left(z, \frac{b_0(1-b_0)r_0}{16\diam\Omega}r\right)}$.
Here,  any pair of point in $\Omega\setminus\overline{B(z, r)}$  can be joined in $\Omega\setminus\overline {B\left(z, \frac{b_0(1-b_0)r_0}{16\diam\Omega}r\right)}$.

Case 2. If $z\in B\left(x_0, \frac{r_0}{8\diam\Omega}r\right)$,
for any $ x\in\Omega\setminus\overline{B(z, r)}$,
$$r-\frac{r_0}{8\diam\Omega}r\leq|x-z|-|x_0-z|\leq|x-x_0|\leq \diam \Omega, $$
so $$r\leq\frac{\diam\Omega}{1-\frac{r_0}{8\diam\Omega}}\leq 2\diam \Omega.$$
Then $$B\left(z, \frac{r_0}{8\diam\Omega}r\right)\subset B\left(x_0, \frac{r_0}{4\diam\Omega}r\right)\subset B\left(x_0, \frac{r_0}{2}\right)\subset B(x_0, r_0)\subset\Omega$$
Similar to the process of proving $b_{x, z, r}>0$ in Lemma \ref{ls48},  we have $\Omega\setminus\overline{B\left(z, \frac{r_0}{8\diam\Omega}r\right)}$ is a connected set. And by
$$\Omega\setminus\overline{B(z, r)}\subset\Omega\setminus\overline{B\left(z, \frac{r_0}{8\diam\Omega}r\right)}, $$
we know any pair of point in $\Omega\setminus\overline{B(z, r)}$ can be joined in $\Omega\setminus\overline{B\left(z, \frac{r_0}{8\diam\Omega}r\right)}$.

Combining above cases,
we get the desired result with $b=\min\left\{\frac{r_0}{16\diam\Omega}, \frac{b_0r_0}{64\diam\Omega}, \frac{b_0(1-b_0)r_0}{16\diam\Omega}  \right\}$.
\end{proof}

\begin{proof}[Proof of Theorem \ref{th1}(ii)]

Let $\Omega\subset\mathbb{R}^n$ be a simply connected planar domain,  or a bounded domain that is quasiconformally equivalent to some uniform domain when $n\leq3$. Assume $\Omega$ supports the $(\phi_\frac{n}{s}, \phi)$-Poincar\'e inequality.

By \cite{bk95,  bsk96}, $\Omega$ has a separation property with $x_0\in\Omega$ and some constant $C_0\geq1$, that is $\forall x\in \Omega$, $\exists$ a curve $\gamma:[0, 1]\to\Omega$, with $\gamma(0)=x, \gamma(1)=x_0$ , and $\forall t\in[0, 1] $, either $\gamma([0, 1])\subset\overline B:=\overline{B(\gamma(t), C_0d(\gamma(t), \Omega^{\complement}))}$, or $\forall y \in \gamma([0, 1])\setminus\overline B$ belongs to  the different component of $\Omega\setminus\overline B$ .
For any $x\in \Omega$ , let $\gamma$ be a curve as above.
By the arguments in \cite{m79}, It suffices to prove
there exists a constant $C>0$ so that
\begin{equation}\label{1.2343}
d(\gamma(t), \Omega^{\complement})\geq C\diam~\gamma([0, t]), ~\forall t\in [0, 1].
\end{equation}
Indeed, \eqref{1.2343} could modify $\gamma$ to get a John curve for $x$.

By Lemma \ref{le421}, $\Omega$ has the LLC(2) property. Let $a=2+\frac{C_0}{b}$, where $b$ is the constant in Lemma \ref{le421}.

For $t\in[0, 1]$.
(1) If $d(\gamma(t), \Omega^{\complement})\geq\frac{d(x_0, \Omega^{\complement})}{a}$, then
$$\gamma([0, t])\subset\Omega\subset B\left(\gamma(t), \frac{ad(\gamma(t), \Omega^{\complement})}{d(x_0, \Omega^{\complement})}\diam \Omega\right).$$
So $$\diam~\gamma([0, t])\leq\frac{2ad(\gamma(t), \Omega^{\complement})}{d(x_0, \Omega^{\complement})}\diam \Omega.$$
and
$$d(\gamma(t), \Omega^{\complement})\geq\frac{d(x_0, \Omega^{\complement})}{2a\diam\Omega}\diam~\gamma([0, t]).$$

(2) If $d(\gamma(t), \Omega^{\complement})<\frac{d(x_0, \Omega^{\complement})}{a}$,  we prove that
$$\gamma([0, t])\subset\overline{B\left(\gamma(t), (a-1)d(\gamma(t), \Omega^{\complement})\right)}.$$
Otherwise,  there exists
$y\in\gamma([0, t])\setminus\overline{B\left(\gamma(t), (a-1)d(\gamma(t), \Omega^{\complement})\right)}$.
By
$$|x_0-\gamma(t)|\geq d(x_0, \Omega^{\complement})-d(\gamma(t), \Omega^{\complement})>(a-1)d(\gamma(t), \Omega^{\complement}), $$
we know $x_0, y\in\Omega\setminus\overline{B\left(\gamma(t), (a-1)d(\gamma(t), \Omega^{\complement}\right)}$,
by Lemma \ref{le421}, $x_0$ and $y$ are contained in the same complement of $\Omega\setminus\overline{B\left(\gamma(t), b(a-1)d(\gamma(t), \Omega^{\complement}\right)}$.
Since $b(a-1)\geq C_0$, then $x_0$ and $y$ are contained in the same complement of $\Omega\setminus\overline{B\left(\gamma(t), C_0d(\gamma(t), \Omega^{\complement}\right)}$, which is in contradiction with the separation property.
Hence
$$\gamma([0, t])\subset\overline{B\left(\gamma(t), (a-1)d(\gamma(t), \Omega^{\complement})\right)}, $$
then $$\diam~\gamma([0, t])\leq2(a-1)d(\gamma(t), \Omega^{\complement}).$$
So $$d(\gamma(t), \Omega^{\complement})\geq\frac{1}{2(a-1)}\diam~\gamma([0, t]).$$

Let $C=\min\left\{\frac{d(x_0, \Omega^{\complement})}{2a\diam\Omega}, \frac{1}{2(a-1)}\right\}$, then \eqref{1.2343}  holds. The proof is completed.

\end{proof}

\medskip
 \noindent {\bf Acknowledgment}. The authors would like to thank  Professor Yuan Zhou for several valuable discussions of this paper.
The authors are partially supported by National Natural Science Foundation of China (No. 12201238) and GuangDong Basic and Applied Basic Research Foundation (Grant No. 2022A1515111056).

\medskip




\begin{thebibliography}{30}

\bibitem{a75}
R.~A. Adams.
\newblock {\em Sobolev spaces}.
\newblock Academic Press, 1975.

\bibitem{refa135}
R.~A. Adams.
\newblock On the {Orlicz-Sobolev} imbedding theorem.
\newblock {\em Journal of Functional Analysis}, 24(3):241--257, 1977.

\bibitem{c97}
A.~Alberico.
\newblock Boundedness of solutions to anisotropic variational problems.
\newblock {\em Communications in Partial Differential Equations},
  36(3):470--486, 2010.

\bibitem{refc3}
A.~Alberico, A.~Cianchi, L.~Pick, and L.~Slav\'ikov\'a.
\newblock Fractional {Orlicz}-{Sobolev} embeddings.
\newblock {\em Journal de Math\'ematiques Pures et Appliqu\'ees}, 149:216--253,
  2021.

\bibitem{b89}
B.~Bojarski.
\newblock Remarks on {Sobolev} imbedding inequalities.
\newblock In {\em Complex {Analysis} {Joensuu} 1987}, pages 52--68. Springer,
  1988.

\bibitem{b82}
J.~Boman.
\newblock {\em $L^p$-estimates for Very Strongly Elliptic Systems}.
\newblock Stockholms Universitet. Matematiska Institutionen, 1982.

\bibitem{bk95}
S.~Buckley and P.~Koskela.
\newblock {Sobolev-Poincar}{\'e} implies {John}.
\newblock {\em Mathematical Research Letters}, 2(5):577--593, 1995.

\bibitem{bsk96}
S.~M. Buckley and P.~Koskela.
\newblock Criteria for imbeddings of {Sobolev-Poincar}{\'e} type.
\newblock {\em International Mathematics Research Notices (IMRN)},
  1996(18):881--901, 1996.

\bibitem{bkl95}
S.~M. Buckley, P.~Koskela, and G.~Lu.
\newblock Boman equals {John}.
\newblock In {\em Proceedings of the XVI Rolf Nevanlinna Colloquium}, pages
  91--99. W. de Gruyter, 1996.

\bibitem{c61}
A.~Calder{\'o}n.
\newblock Lebesgue spaces of differentiable functions and distributions.
\newblock In {\em Proc. Sympos. Pure Math}, volume~4, pages 33--49, 1961.

\bibitem{c96}
A.~Cianchi.
\newblock A sharp embedding theorem for {Orlicz}-{Sobolev} spaces.
\newblock {\em Indiana University Mathematics Journal}, pages 39--65, 1996.
\newblock Publisher: JSTOR.

\bibitem{refc2}
A.~Cianchi.
\newblock A fully anisotropic {Sobolev} inequality.
\newblock {\em Pacific Journal of Mathematics}, 196(2):283--294, 2000.

\bibitem{refb}
A.~Cianchi.
\newblock Optimal {Orlicz}-{Sobolev} embeddings.
\newblock {\em Revista Matem\'atica Iberoamericana}, 20(2):427--474, 2004.

\bibitem{refc1}
A.~Cianchi.
\newblock Higher-order {Sobolev} and {Poincar\'e} inequalities in {Orlicz}
  spaces.
\newblock {\em Forum Mathematicum}, 18(5):745--767, 2006.

\bibitem{refc5}
A.~Cianchi and V.~Maz'ya.
\newblock {Sobolev} inequalities in arbitrary domains.
\newblock {\em Advances in Mathematics}, 293:644--696, 2016.
\newblock Publisher: Elsevier.

\bibitem{cmp11}
D.~V. Cruz-Uribe, J.~M. Martell, and C.~P\'erez.
\newblock {\em Weights, extrapolation and the theory of {Rubio} de {Francia}},
  volume 215.
\newblock Springer Science \& Business Media, 2011.

\bibitem{d16}
B.~Dyda, L.~Ihnatsyeva, and A.~V. V\"{a}h\"{a}kangas.
\newblock On improved fractional {Sobolev-Poincar\'e} inequalities.
\newblock {\em Arkiv f\"{o}r Matematik}, 54(2):437--454, 2016.
\newblock Publisher: Springer.

\bibitem{refe1}
D.~E. Edmunds, R.~Kerman, and L.~Pick.
\newblock Optimal {Sobolev} imbeddings involving rearrangement-invariant
  quasinorms.
\newblock {\em Journal of Functional Analysis}, 170(2):307--355, 2000.
\newblock Publisher: Elsevier.

\bibitem{g88}
D.~Gallardo.
\newblock Orlicz spaces for which the {Hardy-Littlewood} maximal operator is
  bounded.
\newblock {\em Publicacions Matematiques}, pages 261--266, 1988.
\newblock Publisher: JSTOR.

\bibitem{gt}
D.~Gilbarg and N.~S. Trudinger.
\newblock {\em Elliptic partial differential equations of second order}, volume
  224.
\newblock Springer, 1977.

\bibitem{gkz13}
A.~Gogatishvili, P.~Koskela, and Y.~Zhou.
\newblock Characterizations of {Besov} and {Triebel}-{Lizorkin} {Spaces} on
  {Metric} {Measure} {Spaces}.
\newblock {\em Forum Mathematicum}, 25(4):787--819, 2013.

\bibitem{refM1}
L.~Grafakos.
\newblock {\em Classical {Fourier Analysis}}, volume~2.
\newblock Springer, 2008.

\bibitem{refM2}
L.~Grafakos.
\newblock {\em Modern {Fourier Analysis}}, volume 250.
\newblock Springer, 2009.

\bibitem{h01}
P.~Haj\l{}asz.
\newblock {Sobolev} inequalities, truncation method, and {John} domains.
\newblock {\em Rep. Univ. Jyv\"{a}skyl\"{a} Dep. Math. Stat.}, 83:1096--126,
  2001.

\bibitem{refh1}
T.~Heikkinen and H.~Tuominen.
\newblock {Orlicz-Sobolev} extensions and measure density condition.
\newblock {\em Journal of Mathematical Analysis and Applications},
  368(2):508--524, 2010.

\bibitem{h13}
R.~Hurri-Syrj{\"a}nen and A.~V. V{\"a}h{\"a}kangas.
\newblock On fractional {Poincar}{\'e} inequalities.
\newblock {\em Journal d'Analyse Math{\'e}matique}, 120(1):85--104, 2013.

\bibitem{j80}
P.~W. Jones.
\newblock Extension theorems for {BMO}.
\newblock {\em Indiana University Mathematics Journal}, 29(1):41--66, 1980.

\bibitem{j81}
P.~W. Jones.
\newblock Quasiconformal mappings and extendability of functions in {Sobolev}
  spaces.
\newblock {\em Acta Mathematica}, 147:71--88, 1981.

\bibitem{jw78}
A.~Jonsson and H.~Wallin.
\newblock A {Whitney} extension theorem in {$ L^p $} and {Besov} spaces.
\newblock In {\em Annales de l'institut Fourier}, volume~28, pages 139--192,
  1978.

\bibitem{jw84}
A.~Jonsson and H.~Wallin.
\newblock Function spaces on subsets of $\mathbb {R}^n$.
\newblock {\em Mathematical Reports}, 2(1):221, 1984.

\bibitem{refk1}
R.~Kerman and L.~Pick.
\newblock Optimal {Sobolev} imbeddings.
\newblock {\em Forum Mathematicum}, 18(4):535--570, 2006.

\bibitem{refM3}
G.~Leoni.
\newblock {\em A first course in {Sobolev} spaces}.
\newblock American Mathematical Society, 2017.

\bibitem{reflt2}
T.~Liang.
\newblock Fractional {Orlicz--Sobolev} extension/imbedding on {Ahlfors} $ n
  $-regular domains.
\newblock {\em Zeitschrift f{\"u}r Analysis und ihre Anwendungen},
  39(3):245--275, 2020.

\bibitem{reflt3}
T.~Liang.
\newblock A fractional {Orlicz-Sobolev} {Poincar{\'e}} inequality in {John}
  domains.
\newblock {\em Acta Mathematica Sinica}, 37(6):854--872, 2021.

\bibitem{reflt1}
T.~Liang and Y.~Zhou.
\newblock {Orlicz-Besov} extension and {Ahlfors} $ n $-regular domains.
\newblock {\em SCIENTIA SINICA Mathematica}, 51(12):1993--2012, 2020.

\bibitem{m79}
O.~Martio.
\newblock Injectivity theorems in plane and space.
\newblock {\em Ann. Acad. Sci. Fenn. Math. Diss. Ser. A}, 14:384--401, 1978.

\bibitem{m88}
O.~Martio.
\newblock {John} domains, bi-{Lipschitz} balls and {Poincar}{\'e} inequality.
\newblock {\em Rev. Roumaine Math. Pures Appl}, 33(1-2):107--112, 1988.

\bibitem{r83}
Y.~G. Reshetnyak.
\newblock Integral representations of differentiable functions in domains with
  a nonsmooth boundary.
\newblock {\em Sibirskii Matematicheskii Zhurnal}, 21(6):108--116, 1980.

\bibitem{s70}
E.~M. Stein.
\newblock {\em Singular integrals and differentiability properties of
  functions}, volume~2.
\newblock Princeton university press, 1970.

\bibitem{refsh1}
H.~Sun.
\newblock An {Orlicz-Besov Poincar}{\'e} inequality via {John} domains.
\newblock {\em Journal of Function Spaces}, 2019:5234507, 2019.

\bibitem{z123}
Y.~Zhou.
\newblock Criteria for optimal global integrability of {Ha}j{\l}asz-{Sobolev}
  functions.
\newblock {\em Illinois Journal of Mathematics}, 55(3):1083--1103, 2011.

\bibitem{z11}
Y.~Zhou.
\newblock {Haj}{\l}asz-{Sobolev} imbedding and extension.
\newblock {\em Journal of mathematical analysis and applications},
  382(2):577--593, 2011.


\end{thebibliography}
\end{document}